\documentclass[11pt,reqno]{amsart}

\usepackage{xstring}

\newcommand{\globalscale}{0.5}

\usepackage{amsmath,amsfonts,amssymb, amsthm}
\usepackage{bbm}
\usepackage{mathtools}
\usepackage{setspace}
\usepackage{algorithmic}
\usepackage{algorithm}
\usepackage{array}
\usepackage[caption=false,font=normalsize,labelfont=sf,textfont=sf]{subfig}
\usepackage{textcomp}
\usepackage{stfloats}
\usepackage{url}
\usepackage{verbatim}
\usepackage{graphicx}
\usepackage{cite}
\usepackage{enumerate}
\usepackage{changepage}

\usepackage{tikz}
\newcommand\copyrighttext{%
  \footnotesize \textcopyright 2026 IEEE. Personal use of this material is permitted.
  Permission from IEEE must be obtained for all other uses, in any current or future
  media, including reprinting/republishing this material for advertising or promotional
  purposes, creating new collective works, for resale or redistribution to servers or
  lists, or reuse of any copyrighted component of this work in other works.}

\usepackage{caption}
\usepackage[font=small,labelfont=bf]{caption}

\newcommand{\Rset}{\mathbb R}  %ams bold
\newcommand{\Zset}{\mathbb Z}
\newcommand{\Sset}{\mathbb S}

\newcommand{\vect}[1]{\mathrm{vec}(#1)}
\newcommand{\vecti}[1]{\mathrm{vec}^{-1}(#1)}
\newcommand{\svec}[1]{\mathrm{svec}(#1)}
\newcommand{\sveci}[1]{\mathrm{svec}^{-1}(#1)}
\newcommand\norm[1]{\lVert#1\rVert}
\newcommand\Tr[1]{\mathrm{Tr}(#1)}

\newcommand\eigmax[1]{\lambda_{\mathrm{max}}\left(#1\right)}
\newcommand\eigmin[1]{\lambda_{\mathrm{min}}\left(#1\right)}

\theoremstyle{plain}
\newtheorem{thm}{Theorem}
\newtheorem{cor}{Corollary}
\newtheorem{lem}{Lemma}
\newtheorem{assum}{Assumption}
\newtheorem{prop}{Proposition}
\theoremstyle{definition}
\newtheorem{defn}{Definition}
\newtheorem{prd}{Procedure}

\newtheorem{rem}{Remark}

\newcount\Comments
\usepackage{xcolor}
\definecolor{darkgreen}{rgb}{0,0.5,0}
\definecolor{purple}{rgb}{1,0,1}

\newcommand{\keepmarked}[2]{%
  \ifnum\Comments=1
    \textcolor{#1}{#2}%
  \else
    #2%
  \fi
}

\newcommand{\dropmarked}[2]{%
  \ifnum\Comments=1
    \textcolor{#1}{#2}%
  \fi
}

\begin{document}
\title[A Fully Data-Driven Value Iteration for Stochastic LQR]{A Fully Data-Driven Value Iteration for Stochastic LQR: Convergence, Robustness, and Stability}

\author[Leilei Cui, Zhong-Ping Jiang, Petter N.~Kolm, Gregoire G. Macqueron]{
Leilei Cui$^{a}$,  Zhong-Ping Jiang$^{b}$, Petter N.~Kolm$^{c}$, Gr\'egoire G.~Macqueron$^{c,*}$}
\thanks{$^{a}$Department of Mechanical Engineering, University of New Mexico, Albuquerque, New Mexico 87110 USA (e-mail: lcui@unm.edu).}
\thanks{$^{b}$Department of Electrical and Computer Engineering, Tandon School of Engineering, New York University, Brooklyn,
NY 11201 USA (e-mail: zjiang@nyu.edu).}
\thanks{$^{c}$Department of Mathematics, Courant Institute of Mathematical Sciences, New York University, Manhattan, NY 10012 USA (e-mail: petter.kolm@nyu.edu; g.macqueron@nyu.edu).}
\thanks{$^{*}$Corresponding author (g.macqueron@nyu.edu).}
\thanks{\fbox{\parbox{\dimexpr\linewidth-2\fboxsep-2\fboxrule\relax}{\copyrighttext}}}

\def\abstractidth{-1cm}
\begin{adjustwidth}{\abstractidth}{\abstractidth}
\begin{abstract}
Unlike traditional model-based reinforcement learning approaches that estimate system parameters from data, non-model-based data-driven control learns the optimal policy directly from input-state data without any intermediate model identification. Although this direct reinforcement learning approach offers increased adaptability and resilience to model misspecification, its reliance on raw data leaves it vulnerable to system noise and disturbances that may undermine convergence, robustness, and stability. In this article, we establish the convergence, robustness, and stability of value iteration (VI) for data-driven control of stochastic linear quadratic (LQ) systems in discrete-time with entirely unknown dynamics and cost. Our contributions are three-fold. First, we prove that VI is globally exponentially stable for any positive semidefinite initial value matrix in noise-free settings, thereby significantly relaxing restrictive assumptions on initial value functions in existing literature. Second, we extend our analysis to settings with external disturbances, proving that VI maintains small-disturbance input-to-state stability (ISS) and converges within a small neighborhood of the optimal solution when disturbances are sufficiently small. Third, we propose a new non-model-based robust adaptive dynamic programming (ADP) algorithm for adaptive optimal controller design, which, unlike existing procedures, requires no prior knowledge of an initial admissible control policy. Numerical experiments on a ``data center cooling'' problem demonstrate the convergence and stability of the algorithm compared to established methods, highlighting its robustness and adaptability for data-driven control in noisy environments. Finally, we apply the method to dynamic portfolio allocation, demonstrating its practical relevance outside traditional control tasks.
\end{abstract}

\end{adjustwidth}
%
% \copyrightnotice
% {\IEEEpubid{%
% \makebox[\textwidth][c]{%
% \parbox{\textwidth}{\centering\small
% 2162-237X \copyright\ 2026 IEEE. All rights reserved, including rights for text and data mining, and training of artificial intelligence and similar technologies. Personal use is permitted, but republication/redistribution requires IEEE permission. See https://www.ieee.org/publications/rights/index.html for more information.
% }}%
% }
\maketitle
\onehalfspacing
\section{Introduction}
Machine learning (ML) and learning-based methods have become integral to science, engineering, and technology, significantly impacting control theory while introducing new challenges and open problems in the field. Unlike conventional model-based control approaches, which rely on mathematical models derived from first principles, learning-based control learns directly from data and optimizes control strategies. This enables learning-based control to handle complex systems with flexibility and adaptability beyond what manual tuning and traditional model-based strategies offer, thereby enhancing its practical applicability and value. Researchers commonly employ reinforcement learning (RL) techniques to develop learning-based control methods, allowing agents to interact with their environments, make decisions, and minimize cumulative costs \cite{montague1999reinforcement}.

Recent research examines RL algorithms within settings characterized by continuous states and actions, employing the linear quadratic regulator (LQR) as a standard benchmark to evaluate the convergence and robustness of these algorithms \cite{abbasi2011regret,dean2020sample,abbasi2019model,ha2023automated}. Since Kalman introduced the LQR in 1960 \cite{Kalman1960}, subsequent work has demonstrated its theoretical tractability and broad applicability in the controller design of real-world engineering systems \cite{stevens2015aircraft, Huang2022, Cui2023CDC}.

Commonly, RL methods for LQR build upon policy iteration (PI) and value iteration (VI).  PI-based approaches require a stabilizing control policy at the beginning of the learning process \cite{beard1997galerkin,kleinman1969optimal,leake1967construction}. However, this requirement is restrictive in practice. When the system model is not fully accessible, identifying an initial stabilizing controller becomes a significant challenge. Additionally, each PI iteration requires solving a matrix equation, making the method computationally expensive in higher dimensions \cite{howard1960dynamic,kleinman1968riccati,Barto1994,jiang2012computational,bian2014adaptive}. In contrast, VI does not require an initial stabilizing control policy or the solution of a matrix equation at each step. Instead, VI starts from an arbitrary positive definite matrix that serves as an initial approximation of the value function \cite{bellman1966dynamic,bian2016value,Bian2021VI}. In a noise-free environment, and under the assumption of detectability and stabilizability, VI is guaranteed to converge to a stabilizing and optimal controller, provided that each iteration of VI is computed exactly \cite{bertsekas2011dynamic, bian2016value}. These advantages make VI one of the most widely used and thoroughly understood algorithms for solving discounted Markov decision processes (MDPs) \cite{puterman2014markov}. Consequently, VI finds broad application in various engineering domains, including wheel-legged robots \cite{Cui2021robot}, sensory-motor control \cite{Cui2023CDC}, and autonomous vehicles \cite{Huang2022}.

Beyond the idealized case of noise-free environments, a practical learning-based control algorithm must also be robust to perturbations arising from noisy measurements, real-world data, and numerical errors. In data-driven control scenarios, where exact system dynamics are unknown, approximating the value function using sampled data introduces errors due to noise or insufficient data. Consequently, it is essential to determine whether the VI method remains robust to such errors in the learning process. In particular, when subjected to noise and errors, does VI converge to a near-optimal solution? While VI has been extensively studied, its convergence and robustness under these conditions remain unresolved.

Recent work on exact VI establishes exponential convergence only under restrictive assumptions on the initial value function. In particular, \cite{Lee2022TAC} proves exponential convergence when the initial value lies in a semidefinite cone around the optimum, and \cite{song2024convergence} derives local exponential convergence in a neighborhood of the optimum. Relatedly, \cite{lai2025value} studies robustness of inexact VI using an invariant metric and proves convergence to a neighborhood of the optimum under sufficiently small disturbances; however, the result is still local, with bounds that depend on the initial value matrix. 
Finally, several model-free procedures guarantee only near-optimal convergence and may rely on partial model information. For example, \cite{shen2024data} establishes near-optimal convergence for a 
specific LQR class (fast-sampling singular perturbed systems) using partial system knowledge, whereas \cite{fan2024value} considers fully unknown systems but proves only convergence in probability and, as stated in Remark~5 in their article, requires the entire sequence of value iterates to be uniformly bounded, an assumption not guaranteed by the algorithm itself. The authors of \cite{yi2021ALQG,yi2021reinforcementlearning} study data-driven VI algorithms for systems subject to unmeasurable stochastic noise.

In this article, we study the convergence, robustness, and stability of VI for the stochastic LQR problem. First, we show that under noise-free conditions, VI exhibits global exponential convergence to the optimal solution from any positive semidefinite initial value function. Second, in the presence of noise in the VI procedure, we prove that VI converges to a small neighborhood of the optimal solution, provided that the noise is sufficiently small. Third, as a direct application of these results to adaptive dynamic programming (ADP) \cite{jiang2020learning,book_Jiang}, we develop a non-model-based robust least-squares value iteration (R-LSVI) algorithm for adaptive optimal controller design. We empirically demonstrate its convergence, stability, and adaptability, and compare it with several established methods on a stochastic LQR problem -- a stylized ``data center cooling'' scenario frequently used to benchmark reinforcement learning approaches \cite{gao2014machine, tu2018least, abbasi2019model, dean2020sample}. Finally, we demonstrate the effectiveness of our approach on a dynamic portfolio allocation problem \cite{garleanu2013dynamic}, a real-world problem from finance.

We emphasize that the convergence and robustness results presented in this article significantly relax the restrictive conditions on the initial value function and noise that are typically assumed in the existing literature \cite{song2024convergence,Lee2022TAC,zhang2023revisiting,lai2025value}. A central technique in our analysis is the adoption of input-to-state stability (ISS) from modern control theory \cite{Sontag1989,Sontag2008}, specifically focusing on its recent variant, ``small-disturbance ISS'' \cite{pang2021robust,Cui2024SCL}. This framework facilitates a comprehensive analysis of VI’s convergence and robustness in the presence of disturbances and errors.

The remainder of this article is organized as follows. In Section \ref{sec:preliminaries-and-problem-formulation}, we present the necessary preliminaries, formulate the problem setup, and summarize our main results. Section \ref{sec:robust-least-squares-VI} introduces the R-LSVI algorithm, a fully data-driven method for constructing estimates of the value function directly from input-state data. In Section \ref{sec:convergence-analysis}, we present a comprehensive convergence analysis of VI. We establish the global exponential stability of the exact VI method. This result is key to  demonstrating that inexact VI achieves small-disturbance ISS. The convergence of the R-LSVI algorithm follows from it being a special case of inexact VI. In Section \ref{sec:empirical-results}, we numerically demonstrate the convergence, stability, and adaptability of R-LSVI, comparing it to existing methods using a stochastic LQR problem inspired by a ``data center cooling'' scenario. In Section \ref{sec:portfolio-allocation}, we apply R-LSVI to a dynamic portfolio allocation problem, a practical application in finance. Section \ref{sec:conclusions} concludes. We relegate some technical proofs to the appendix. 

\noindent\textit{Notation:} $\Rset$ (resp.~$ \Rset_{+}$) denotes the set of (nonnegative) real numbers; $\Zset_{+}$ is the set of nonnegative integers; $\mathbb{S}^{n}$ is the set of real symmetric matrices of dimension $n$; $\Sset^{n}_{+}$ (resp. $\Sset^{n}_{++}$) is the set of real symmetric and positive semi-definite (definite) matrices of dimension $n$; $\otimes$ denotes the Kronecker product; $I_{n}$ denotes the identity matrix of dimension $n$ while $0_{m \times n}$ denotes the $m \times n$ zero matrix; $\|\cdot\|_{2}$ denotes the Euclidean norm for vectors and the spectral norm for matrices; for a function $u: \mathbb{F} \rightarrow \Rset^{n \times m}$, $\|u\|_{\infty}$ denotes its $l^{\infty}$-norm when $\mathbb{F}=\Zset_{+}$, and $L^{\infty}$-norm when $\mathbb{F}=\Rset$. $\mathbf{1}_{[0,T]}$ denotes an indicator function, that is $\mathbf{1}_{[0,T]}(t)=1$ if $t \in [0,T]$, and $\mathbf{1}_{[0,T]}(t)=0$ otherwise. For matrices $X \in \Rset^{m \times n}, Y \in \mathbb{S}^{m}$, define $\vect{X} \triangleq \left[X_{1}^{\top}, X_{2}^{\top}, \ldots, X_{n}^{\top}\right]^{\top}\in\Rset^{nm}$, $\svec{Y} \triangleq \left[y_{11}, \sqrt{2} y_{12}, \ldots, \sqrt{2} y_{1 m}, y_{22}, \sqrt{2} y_{23}, \ldots, \sqrt{2} y_{m-1, m}, y_{m, m}\right]^{\top}$ $\in\Rset^{\frac{1}{2} m(m+1)}$, where $X_{i}$ is the $i^{th}$ column of $X$. $\vecti{\cdot}$ and $\sveci{\cdot}$ are the operations such that $\vecti{\vect{X}} = X$, and $\sveci{\svec{Y}} = Y$. For $n\in\Zset_{+}$, define $ \tilde{n} \triangleq \frac{n(n+1)}{2}$. For a vector $v \in \Rset^{n}$, define $\tilde{v} \triangleq \svec{vv^{\top}} \in \Rset^{\tilde{n}}$. For a vector $v\in\Rset^{n}$, define $\mathcal{E}(v) \triangleq \left[v_1,\dots,v_{n-1}\right]^{\top}\in\Rset^{n-1}$ as the vector obtained by deleting the last element of $v$. For $Y\in\mathbb{S}^{n}$,
$\eigmin{Y}$ ($\eigmax{Y}$) denotes its smallest (largest) 
algebraic eigenvalue. For $Y_{1}, Y_{2}\in \mathbb{S}^{n}$, $Y_{1} \preceq Y_{2}$ means that $Y_{2} - Y_{1} \in\Sset^{n}_{+}$. For $S\in \Sset^{n}_{+}$, $S^{\frac{1}{2}}$ denotes the upper triangular Cholesky factor of $S$. For a matrix $Z\in \Rset^{n \times n}$, $Z^{\dagger}$ denotes its Moore-Penrose pseudoinverse. For $Z \in \Rset^{m \times n}$, define $\mathcal{B}_{r}(Z) \triangleq \left\{X \in \Rset^{m \times n}, \norm{X-Z}_{2}<r\right\}$ to be the open ball of radius $r$ centered at $Z$, and $\overline{\mathcal{B}}_{r}(Z)$ its closure, similarly $\mathcal{B}^{\infty}_{r}(Z)$ denotes the ball with respect to the infinity norm.

\section{Preliminaries and Problem Formulation}\label{sec:preliminaries-and-problem-formulation}

\subsection{Exponential and Input-to-State Stability}
\noindent In this article, we use exponential stability and input-to-state stability to analyze the convergence of VI. Below we recall the key definitions (see  \cite{Sontag2008,Cui2024SCL} for details). \\
\indent We consider the nonlinear system 
\begin{equation}
\label{eq:dynamic_gen}
P_{i+1}=f(P_{i}, \Delta_{i}) 
\end{equation}
where $P_{i} \in \Rset^{n}$ is the system state, $\Delta_{i} \in \Rset^{m}$ is the perturbation, $f: \Rset^{n} \times \Rset^{m} \rightarrow \Rset^{n}$ is continuous, and $P^{*}$ is an equilibrium of $P_{i+1}=f\left(P_{i}, 0\right)$ when $\Delta_{i}=0$ for all $i \in \Zset_{+}$.
\begin{defn}
For the system \eqref{eq:dynamic_gen} with $\Delta_{i}=0$ for all $i \in \Zset_{+}$, an equilibrium $P^{*}$ is locally exponentially stable if there exists a $\delta>0$, such that for some $a>0$ and $0<b<1$,
\begin{equation}
    \norm{P_{i}-P^{*}}_{2} \leq a b^{i}\norm{P_{0}-P^{*}}_{2}
\end{equation}
for all $P_{0} \in \mathcal{B}_{\delta}\left(P^{*}\right)$ and all $i \in \mathbb{Z}_+$. If $\delta=\infty$, then $P^{*}$ is a globally exponentially stable equilibrium.
\end{defn}
\begin{defn}
    A function $\gamma: \Rset_{+} \rightarrow \Rset_{+}$ is said to be of class $\mathcal{K}$ if it is continuous, strictly increasing and vanishes at the origin. A function $\beta: \Rset_{+} \times \Zset_{+} \rightarrow \Rset_{+}$ is said to be of class $\mathcal{KL}$ if $\beta(\cdot, i)$ is of class $\mathcal{K}$ for every fixed $i >0$, and for every fixed $r \geq 0$, $\beta(r, i)$ decreases to $0$ as $i \rightarrow \infty$.
\end{defn}
\begin{defn}
    The system \eqref{eq:dynamic_gen} is small-disturbance input-to-state stable (ISS) if there exist some $d>0$, $\beta \in \mathcal{KL}$ and $\gamma \in \mathcal{K}$, such that for all $\Delta$ bounded by $d$ (i.e.~$\left\|\Delta\right\|_{\infty}<d$) and all initial states $P_{0}\in\Rset^{n}$, $P_{i}$ satisfies 
    \begin{equation}
        \norm{P_{i}-P^{*}}_{2} \leq \beta\left( \norm{P_{0}-P^{*}}_{2}, i \right) +\gamma\left(\norm{\Delta}_{\infty}\right) ,\, \forall i \in \mathbb{Z}_+.
    \end{equation}
\end{defn}
Exponential stability implies not only convergence of system \eqref{eq:dynamic_gen}, but also characterizes its rate of convergence. When the input signal is non-zero, the input-to-state stability describes how the solution of \eqref{eq:dynamic_gen} is affected by the input signal. In particular, small-disturbance ISS ensures that the distance between the system state and the equilibrium remains uniformly bounded whenever the disturbance is sufficiently small, regardless of the initial state.

\subsection{The Stochastic Linear Quadratic Regulator Problem in Discrete-Time}
\noindent We consider a stochastic discrete-time linear system of the form
\begin{equation}
\label{eq:dynamic}
    x_{t+1}=A x_{t}+B u_{t} + C\epsilon_{t+1}
\end{equation}
where $x_{t} \in \Rset^{n}$ is the system state, $u_{t} \in \Rset^{m}$ is the control input, $\epsilon_{t+1} \in \Rset^{p}$ is an additive disturbance, $A \in \Rset^{n \times n}, B \in \Rset^{n \times m}$ and $C \in \Rset^{n \times p}$.  We assume that the initial state $x_{0} \in \Rset^{n}$ is given,  $\mathbb{E}[\epsilon_{t}] = 0_{p}$, $\mathbb{V}[\epsilon_{t}] = I_{p}$, and $\epsilon_{t},\epsilon_{t'}$ are independent for all $t'\neq t$. The stochastic LQR problem seeks a controller $u$ that minimizes the cost 
\begin{equation}
\label{eq:functional}
J\left(u\right)= \lim_{T\to\infty} \frac{1}{T} \mathbb{E}\left[\sum_{t=0}^{T-1}c\left(x_t, u_t\right) \right]
\end{equation}
where $c(x_t,u_t) \triangleq x_{t}^{\top} S x_{t}+u_{t}^{\top}Ru_{t}$ with $S\in\Rset^{n\times n}$ and $R\in\Rset^{m\times m}$. 
\begin{assum}
    \label{assum:hyp_control}
    $(A,B)$ is stabilizable, and $S\in\Sset^{n}_{++}$, $R\in\Sset^{m}_{++}$.
\end{assum}
Under Assumption~\ref{assum:hyp_control}, the stochastic LQR problem admits a unique optimal controller $u^{*}_{t}=-K^{*} x_{t}$ with an optimal cost $J\left(u^{*}\right) = \Tr{C^{\top} P^{*} C}$, where $K^{*}\in \Rset^{m\times n}$ satisfies
\begin{equation}
\label{eq:K_star_formula}
    K^{*}=\left(R+B^{\top} P^{*} B\right)^{-1} B^{\top} P^{*} A
\end{equation}
and $P^{*}\in\Sset^{n}_{++}$ is the unique solution of the discrete algebraic Riccati equation (DARE)
\begin{equation}
\label{eq:riccati-eqn}
    P= A^{\top} P A - A^{\top} P B\left(R+B^{\top} P B\right)^{-1} B^{\top} P A+ S\,.  
\end{equation}
In addition, $A-BK^{*}$ is Schur stable, i.e.~the spectral radius $\rho(A-BK^{*}) < 1$ (see, for example,  \cite[Section 2.4]{lewis2012optimal}). For any given $P\in\Sset_{+}^{n}$, the resulting control is
\begin{equation}
\label{eq:K_formula}
    K_{P}=\left(R+B^{\top} P B\right)^{-1} B^{\top} P A 
\end{equation}
where the subscript in $K_{P}$ indicates its dependence on the value matrix $P$. For convenience, we write $K_{i}$ in place of $ K_{P_{i}}$ throughout the article.

We introduce the Hamiltonian operator $\mathcal{Q}:\Rset^{n\times n}\to \Rset^{l\times l}$,  with $l = n+m$,
\begin{align}
    \label{eq:hamiltonian_Q}
    \mathcal{Q}\left(P\right)&\triangleq \left(\begin{array}{cc}
    A^{\top} PA+S & A^{\top} P B \\
    B^{\top} P A & B^{\top} P B+R
    \end{array}\right)
    =  \left(\begin{array}{cc}
    \left[\mathcal{Q}(P)\right]_{xx} & \left[\mathcal{Q}(P)\right]_{ux}^{\top} \\
    \left[\mathcal{Q}(P)\right]_{ux} & \left[\mathcal{Q}(P)\right]_{uu}\\
    \end{array}\right)
\end{align}
which plays a fundamental role in the stochastic LQR problem. We define by $\mathcal{H}$ the $[\cdot]_{uu}$-Schur complement operator, given by 
\begin{equation}
    \mathcal{H}\left(Q\right) \triangleq \left[{Q}\right]_{xx}  - \left[{Q}\right]_{ux}^{\top} \left[{Q}\right]_{uu}^{-1}\left[{Q}\right]_{u x}\in\Rset^{n\times n}
\end{equation}
for any $Q\in\Rset^{l\times l}$ where $\left[{Q}\right]_{uu}$ is invertible. With this notation, the DARE \eqref{eq:riccati-eqn} becomes
\begin{equation}
    \label{eq:P_formula_Q}
    P = \mathcal{H}\left(\mathcal{Q}(P)\right)\,,
\end{equation}
and the greedy control gain \eqref{eq:K_formula} is
\begin{equation}
\label{eq:K_formula_Q}
    K_{P} = \left[\mathcal{Q}\left(P\right)\right]_{uu}^{-1}\left[\mathcal{Q} \left(P\right)\right]_{u x}\,.
\end{equation}
Thus, $K^{*} = K_{P^{*}}$, where $P^{*}\in\Sset^{n}_{++}$ is the unique solution of \eqref{eq:P_formula_Q}. Since $P\in\Sset_{+}^{n}$ and $R\in\Sset_{++}^{m}$, it follows that $\left[\mathcal{Q}\left(P\right)\right]_{uu}\in\Sset_{++}^{m}$ and is therefore invertible. \\
The formulas \eqref{eq:P_formula_Q}--\eqref{eq:K_formula_Q} capture the core principle behind the derivation of the Riccati equation. By eliminating the control variables via the Schur complement $\mathcal{H}$, we obtain an effective state cost that determines the optimal closed-loop behavior of the system.

To solve the stochastic LQR problem \eqref{eq:dynamic}--\eqref{eq:functional} we must solve the DARE \eqref{eq:riccati-eqn} or equivalently \eqref{eq:P_formula_Q}. When the system dynamics are known, this can be accomplished using Bertsekas’s exact value iteration \cite{bertsekas2011dynamic}, as outlined in Procedure \ref{prd:exact_VI}. 

\begin{prd}[Exact Value Iteration]
\label{prd:exact_VI}
\ \begin{enumerate}
    \item Choose an initial value matrix $P_{0}\in\Sset^{n}_{+}$, and let $i=0$.
    \item Update the  value matrix via
        \begin{equation}
        \label{eq:exact_VI_2}
        P_{i+1}=
        \mathcal{H}\left(Q_{i}\right)
        \end{equation}
    where $Q_{i} \triangleq \mathcal{Q}\left(P_{i}\right)$.
    \item Set $i \leftarrow i+1$ and go back to Step (2).
\end{enumerate}
\end{prd}

We note that at any step of Procedure \ref{prd:exact_VI}, the control gain $K_{i}$ associated with the value matrix $P_{i}$ can be explicitly computed as 
\begin{equation}
K_{i} =\left[\mathcal{Q}\left(P_{i}\right)\right]_{uu}^{-1}\left[\mathcal{Q} \left(P_{i}\right)\right]_{u x}. 
\end{equation}
Unlike Hewer’s policy iteration \cite{hewer1971iterative},  Procedure \ref{prd:exact_VI} converges even if $K_i$ is not stabilizing, as formally stated in the following theorem. This follows directly from \cite[Proposition 4.4.1]{bertsekas2011dynamic}.
\begin{thm}[Convergence of Exact VI]
\label{thm:exact_VI}
Under Assumption~\ref{assum:hyp_control}, it holds for Procedure \ref{prd:exact_VI} that 
\begin{enumerate}
    \item[(i)] $\lim _{i \rightarrow \infty} P_{i}=P^{*}$.
    \item[(ii)]$\lim _{i \rightarrow \infty} K_{i}=K^{*}$.
\end{enumerate}
\end{thm}

\subsection{Problem Formulation}
\noindent Theorem \ref{thm:exact_VI} states that exact VI (see Procedure \ref{prd:exact_VI}) finds a sequence of approximations that are guaranteed to converge to the true solution $P^{*}$ and $K^{*}$. However, implementing exact VI requires evaluating $\mathcal{Q}$ defined in \eqref{eq:hamiltonian_Q}, which requires precise knowledge of $A$, $B$, $S$ and $R$.

In the setting of learning-based control, system dynamics are unknown and must be approximated using data-driven approaches. In some situations, the cost function may also be partially or completely unknown. Moreover, data employed in these approaches can be affected by disturbances such as system noise and numerical errors. Consequently, in practice each VI update yields an approximate value matrix $\hat{P}_{i+1}$. To incorporate errors and noise, we formulate the following inexact version of the VI.

\begin{prd}[Inexact Value Iteration]~\\[-4mm]
    \label{prd:inexact_VI}
    \begin{enumerate}
    \item Choose an initial value matrix $\hat{P}_{0}\in\Sset^{n}_{+}$, and let $i=0$.   
    \item Update the value matrix via
    \begin{align*}
        \hat{P}_{i+1}&=\mathcal{H}(\mathcal{Q}(\hat{P}_{i})) + \Delta_{i}\,.
    \end{align*}
    \item Set $i \leftarrow i+1$ and go back to Step (2).
\end{enumerate}
\end{prd}

In Procedure \ref{prd:inexact_VI}, $\Delta_{i} \in \Sset^{n}$ denotes external disturbances, while $\hat{P}_{i+1}$ is used to distinguish the sequence generated by the inexact value iteration from that of the exact value iteration. We always enforce the condition $\Delta_{i} \in \Sset^{n}$ by symmetrizing $\hat{P}_{i+1} = \frac{1}{2}(\hat{P}_{i+1} + \hat{P}_{i+1}^{\top})$.

\begin{rem}
\label{rem:inexact_vi_framework}
The inexact VI is general in the sense that any VI for stochastic LQR subject to errors or disturbances can be expressed as in Procedure \ref{prd:inexact_VI}, allowing us to analyze convergence, robustness, and stability under a broad range of perturbations. We provide two examples to illustrate this point. \\
First, consider the case where we do not have the exact system and cost matrices, but some approximation $(\hat{A}, \hat{B}, \hat{S}, \hat{R})$ thereof.  Then $\Delta_{i}$ takes the form
\begin{equation}
    \Delta_{i} =  \mathcal{H}(\hat{\mathcal{Q}}(\hat{P}_{i})) - \mathcal{H}({\mathcal{Q}}(\hat{P}_{i})) 
\end{equation}
where
\begin{equation}
    \hat{\mathcal{Q}}\left(P\right)\triangleq \left(\begin{array}{cc}
    \hat{A}^{\top} P\hat{A}+\hat{S} & \hat{A}^{\top} P \hat{B} \\
    \hat{B}^{\top} P \hat{A} & \hat{B}^{\top} P \hat{B}+\hat{R}
    \end{array}\right)\,.
\end{equation}
Second, consider the R-LSVI algorithm (see Algorithm \ref{alg:R-LSVI}, below) where the source of disturbances arises from approximating the matrix $\mathcal{Q}(\hat{P}_{i})$. Then $\Delta_{i}$ takes the form
\begin{equation}
    \Delta_{i} =  \mathcal{H}(\hat{\mathcal{Q}}(\hat{P}_{i})) - \mathcal{H}(\mathcal{Q}(\hat{P}_{i})) 
\end{equation}
where $\hat{\mathcal{Q}}(\cdot)$ is defined in Section \ref{sec:robust-least-squares-VI} as
\begin{equation}
    \hat{\mathcal{Q}}\left(P\right) \triangleq \sveci{\mathcal{E}(\Theta_{T}^{\dagger}\Psi_{T}\svec{P} + \Theta_{T}^{\dagger}\Xi_{T})}\,.
\end{equation}
In both examples, the error arises from approximating the Hamiltonian operator $\mathcal{Q}$, but the nature of the approximation differs significantly between the two. In the first example, known as the nominal control approach, the system and cost matrices $(A, B, S, R)$ are approximated first, and then the Hamiltonian is constructed in a model-based manner. In contrast, in the second example, based on adaptive dynamic programming (ADP), we directly approximate the Hamiltonian, resulting in a non-model-based approach. This distinction affects the structure of the perturbations and the resulting error propagation in the inexact value iteration process. In this article, we establish sufficient conditions on $\{\Delta_{i}\}_{i=0}^{\infty}$ to guarantee the convergence of sequences generated by inexact VI to a neighborhood of the true solution $P^{*}$. 
\end{rem}

\subsection{Main Contributions}
\noindent Our contributions are threefold:\\
\noindent\textbf{Result 1}: We show that the exact VI (see Procedure \ref{prd:exact_VI}) is globally exponentially stable for any $P_{0}\in\Sset^{n}_{+}$.\\
\noindent\textbf{Result 2}: We prove that the inexact VI (see Procedure \ref{prd:inexact_VI}) is small-disturbance ISS.\\
\noindent\textbf{Result 3}: We propose a new non-model-based adaptive dynamic programming algorithm based on VI to determine the optimal controller in \eqref{eq:K_star_formula} using only input-state data, with convergence guarantees.

In addition, we evaluate the proposed algorithm empirically on the stochastic LQR ``data center cooling'' benchmark, which is frequently used to assess reinforcement-learning methods (e.g.~\cite{gao2014machine,tu2018least,abbasi2019model,dean2020sample}), and study its convergence, stability, and adaptability. Our simulations compare the proposed algorithm with several alternatives, including both model-based PI and VI  approaches \cite{recht2019tour,dean2020sample}, a model-free policy gradient method \cite{montague1999reinforcement}, and two ADP policy iteration procedures \cite{krauth2019finite, pang2021robust}. The results highlight the superior adaptability and stability of our approach, even when costs are non-quadratic.

\section{Robust Least-Squares VI}\label{sec:robust-least-squares-VI}
\noindent We introduce the robust least-squares VI (R-LSVI) algorithm in this section, a fully data-driven approach to estimate $\hat{P}_{i}$ in Procedure \ref{prd:inexact_VI} from input-state data. Section \ref{sec:convergence-of-the-robust-VI} provides a rigorous convergence analysis. 

We define the filtration $\{\mathcal{F}_t\}_{t \geq 0}$ by $\mathcal{F}_t = \sigma(x_0, \ldots, \allowbreak x_t)$, where $\sigma(\cdot)$ denotes the $\sigma$-algebra generated by the state sequence up to time $t$.

Because $(A,B,S,R)$ are unknown, the behavior policy $u_t$  we choose for data collection may produce an unstable closed loop and drive the state to infinity. To prevent divergence and numerical ill-conditioning, we reset the state to $0_{n}$ whenever $\norm{x_{t}}_{\infty}$ exceeds $d$. We formalize this mechanism in the following proposition.

\begin{prop}
\label{prop:reset}
Define the reset dynamics with initial condition $x_{0} = 0_{n}$
\begin{align}
\label{eq:dyn_reset1}
X_{t+1} &= Ax_t + Bu_t + C\epsilon_{t+1}, \quad  \epsilon_t \underset{iid}{\sim}\mathcal N(0_p, I_p) \\
\label{eq:dyn_reset2}
x_{t+1} &=  \mathbbm{1}_{{\norm{X_t}_{\infty} \le d}}X_{t+1}, \quad d>0
\end{align}
and consider a behavior policy of the form
\begin{equation}
\label{eq:behavior_policy}
    u_t = -K_c x_t + \eta_{t+1}, \quad \eta_t \underset{iid}{\sim}\mathcal N(0_m, \Sigma_\eta), \quad \Sigma_\eta \succ 0
\end{equation}
where $K_c$ is not necessarily stabilizing. 
Then the Markov chain $(x_t)_{t\geq 0}$ admits a unique invariant probability measure $\pi$.
\end{prop}
The proof relies on the irreducibility and positive recurrence of $(x_t)_{t\ge0}$ (via regeneration at $0_n$), which ensure the existence and uniqueness of the invariant measure $\pi$; see Appendix~\ref{apdx:reset} for details.

We now derive the least squares formulation used in the construction of $\hat{\mathcal{Q}}(P)$. For any P $\in\Sset^{n}_{+}$, recall $X_{t+1} = A x_{t} + B u_{t} + C\epsilon_{t+1}$, and hence
\begin{align}
\mathbb{E}\left[X_{t+1}^{\top}PX_{t+1} + c_{t} | \mathcal{F}_{t}\right] 
&= (A x_{t} + B u_{t})^{\top}P( A x_{t} + B u_{t})   + x_{t}^{\top}Sx_{t} +u_{t}^{\top}Ru_{t}\notag
\\&\quad\qquad+ \mathbb{E}\left[(C\epsilon_{t+1})^{\top}P(C\epsilon_{t+1}) | \mathcal{F}_{t}\right]\\
&=\left(\begin{array}{c}
    x_{t} \\
    u_{t}
    \end{array}\right)^{\top}\mathcal{Q}(P)\left(\begin{array}{c}
    x_{t} \\
    u_{t}
    \end{array}\right) + \mu(P) 
\end{align}
where $\mu(P) = \Tr{C^{\top}PC}$. Setting $y_{t} = [x_{t}^{\top}, u_{t}^{\top}]^{\top}$, we obtain
\begin{equation}
\label{eq:TD_vec}
    \mathbb{E}\left[\tilde{X}_{t+1}^{\top}\svec{P} |\mathcal{F}_{t}\right]  + c_{t}  = \tilde{y}_{t}^{\top}\svec{\mathcal{Q}\left(P\right)} + \mu\left(P\right) \,.
\end{equation}
Defining $z_{t} \triangleq \left[\tilde{y}_{t}^{\top},1\right]^{\top}$, multiplying both sides of \eqref{eq:TD_vec} by $\frac{z_{t}}{\alpha_{t}^{2}}$ with $\alpha_{t} = \norm{z_{t}}_{\infty}$, and taking the expectation with respect to the invariant distribution $\pi$ yields 
\begin{equation}
    \label{eq:TD_VI_LS}
    \Theta\left(\begin{array}{l}
    \svec{\mathcal{Q}\left(P\right)}\\
            ~~~~\mu\left(P\right)
    \end{array}\right) = \Psi\svec{P} + \Xi
\end{equation}
where
\begin{equation}
        \label{eq:theta-psi-ksi}
         \Theta \triangleq \mathbb{E}_{\pi}\left[\frac{z_t z_{t}^{\top}}{\alpha_{t}^{2}}\right], \quad
         \Psi \triangleq \mathbb{E}_{\pi}\left[\frac{z_t \tilde{X}_{t+1}^{\top}}{\alpha_{t}^{2}}\right], \quad \Xi \triangleq \mathbb{E}_{\pi}\left[\frac{z_t c_{t}}{\alpha_{t}^{2}}\right].
\end{equation}

\begin{assum}
\label{assum:invertible1}
There exists $c>0$,  such that for $\phi_t \triangleq \frac{z_t}{\alpha_t}$ we have
\begin{equation}
    \Theta =  \mathbb{E}_\pi[\phi_t\phi_t^\top] \succeq c I_{\tilde{l}+1}\,.
\end{equation}
\end{assum} 
\begin{rem}
Assumption \ref{assum:invertible1} is consistent with the persistent excitation condition widely used in adaptive control \cite{aastrom1995adaptive,willems2005note} and is similarly prevalent in various RL methods (see, for example,  \cite{book_Jiang}, \cite{kamalapurkar2018reinforcement}, \cite{kiumarsi2017optimal}).
\end{rem}

For a given $P$ on the right-hand side of \eqref{eq:TD_VI_LS}, we can estimate $\mathcal{Q}\left(P\right)$ on the left hand side from collected data via least-squares (LS), leading to
\begin{equation}
    \label{eq:TD_VI_LS_approx}
    \left(\begin{array}{l}
    \svec{\hat{\mathcal{Q}}\left(P\right)}\\
            ~~~~~\hat{\mu}\left(P\right)
    \end{array}\right) = \Theta_{T}^{\dagger}\Psi_{T}\svec{P} + \Theta_{T}^{\dagger}\Xi_{T}
\end{equation}
where the data matrices
\begin{equation}
    \label{eq:data_matrices}
    \Theta_{T} \triangleq \frac{1}{T}\sum_{t=1}^{T}\frac{z_t z_{t}^{\top}}{\alpha_{t}^{2}},~ 
    \Psi_{T} \triangleq \frac{1}{T}\sum_{t=1}^{T}\frac{z_t \tilde{X}_{t+1}^{\top}}{\alpha_{t}^{2}},~
    \Xi_{T} \triangleq \frac{1}{T}\sum_{t=0}^{T-1} \frac{z_{t}c_{t}}{\alpha_{t}^{2}}
\end{equation}
are sample estimators of  \eqref{eq:theta-psi-ksi}. As the chain $(x_{t})_{t>0}$ is irreducible and positive recurrent, by the ergodic theorem, the data matrices in \eqref{eq:data_matrices} converge with increasing data size, i.e.
\begin{equation}
\lim _{T \rightarrow \infty} \Theta_{T}=\Theta,\quad
\lim _{T \rightarrow \infty} \Psi_{T}=\Psi,\quad
\lim _{T \rightarrow \infty} \Xi_{T}=\Xi  \,.
\end{equation}
\begin{rem}
The rescaling by $\alpha_t=\norm{z_t}_{\infty}$ in \eqref{eq:theta-psi-ksi} and \eqref{eq:data_matrices} acts as a numerical preconditioning step. Because $\norm{x_t}_{\infty}\in[0,d]$ (with $d$ potentially large), the dynamic range of the data can be wide. Without normalization, large-magnitude samples dominate the averages and yield ill-conditioned estimates of  \eqref{eq:data_matrices}.
\end{rem}
We now present the R-LSVI algorithm in Algorithm \ref{alg:R-LSVI}. The method operates entirely without access to system matrices $(A,B,S,R)$. 
We sample data by applying the behavior policy $u_t$ defined in \eqref{eq:behavior_policy} and observing $(X_{t},x_t,c_t)$ along trajectories of \eqref{eq:functional}, \eqref{eq:dyn_reset1} and \eqref{eq:dyn_reset2}, from which we construct the data matrices \eqref{eq:data_matrices}. The subsequent steps follow Procedure \ref{prd:exact_VI}, with $\hat{Q}_{i}$, computed via \eqref{eq:TD_VI_LS_approx}, replacing  $Q_{i}$. Finally, we obtain the control gain $\hat{K}$ using formula \eqref{eq:K_formula_Q}, as it is the primary quantity of practical interest. Note that, unlike the PI Procedure, we do not need the knowledge of a stabilizing initial controller $K_{0}$.

\begin{algorithm}[h]
\caption{Robust Least-Squares Value Iteration (R-LSVI)}
\label{alg:R-LSVI}
\begin{algorithmic}
    \STATE \textbf{input:} Choose $\hat{P}_{0}\in\Sset^{n}_{+}$, exploration-noise  covariance $\Sigma_{\eta}$, number of iterations $I_{\text{max}}$, time horizon $T$, the reset bound $d>0$, and let $i=0$.
    \STATE Collect data $(X_{t}, x_{t}, u_{t})$ from \eqref{eq:dyn_reset1}, \eqref{eq:dyn_reset2} and \eqref{eq:behavior_policy}.
    \STATE Construct $\Theta_{T}$, $\Psi_{T}$ and $\Xi_{T}$.
    \WHILE{ $i < I_{\text{max}}$} 
        \STATE $\hat{Q}_{i} = \sveci{\mathcal{E}(\Theta_{T}^{\dagger}\Psi_{T}\svec{\hat{P}_{i}} + \Theta_{T}^{\dagger}\Xi_{T})}$
        \vspace{0.5mm}
        \vspace{0.5mm}
        \STATE $\hat{P}_{i+1}=\mathcal{H}(\hat{Q}_{i})$
        \STATE Update $i \leftarrow i+1$.
    \ENDWHILE
    \STATE \textbf{output:} $\hat{K}=[\hat{Q}_{I_{\text{max}}}]_{uu}^{\dagger}[\hat{Q}_{I_{\text{max}}}]_{u x}$
\end{algorithmic}
\end{algorithm}

\begin{rem} 
This procedure significantly sharpens previous results. Specifically, the law of $(x_t)$ converges to the invariant distribution $\pi$,  leading to almost sure convergence of the data matrices, whereas prior work established only convergence in probability (e.g.~\cite{fan2024value}).
\end{rem}
We make the following remarks regarding the data matrices defined in the main article by equation \eqref{eq:data_matrices} in the R-LSVI algorithm:
\begin{enumerate}[(i)]
    \item  The algorithm computes $\Theta_T^\dagger \Psi_T$ and $\Theta_T^\dagger \Xi_T$ once and reuses them in every iteration, which improves computational efficiency. The overall complexity is dominated by computing the pseudoinverse of $\Theta_T$ and the matrix products performed at each step. The overall complexity is $O(T\tilde{l}^2 + I_{\max}\tilde{l}^2)$.
    \item  R-LSVI constructs the data matrices entirely from input-state data and, unlike the approaches in \cite{bian2020model,jiang2014adaptive,jiang2020learning,bian2016adaptive}, does not rely on observing system disturbances $\{\epsilon_{t}\}_{t\geq0}$.
    \item The data matrices are smaller in size than the ones typically found in the ADP literature (see, for example,  \cite{jiang2012computational}, \cite{pang2022reinforcement}). For instance, $\Theta_{T}$ is of size $(\tilde{l} + 1)^{2}$ rather than  $(\tilde{l} + l + 1)^{2}$ . This dimension reduction arises from considering
    $\left(\begin{array}{l}
    \svec{\mathcal{Q}\left(P\right)}\\
            ~~~~\mu\left(P\right)
    \end{array}\right)$ rather than
    $\left(\begin{array}{ll}
    \svec{\mathcal{Q}\left(P\right) } & 0_{(n+m)} \\
            0_{(n+m)}^{\top}&~~~\mu\left(P\right)
    \end{array}\right)$ in the LS problem defined in equation \eqref{eq:TD_VI_LS} of the main article, leading to a mathematically equivalent but computationally more efficient formulation.
\end{enumerate}
% Additional remarks on the algorithm complexity can be found in Appendix~\ref{apdx:rem_complexity}.

\section{Convergence Analysis}
\label{sec:convergence-analysis}
\noindent Detailed proofs of all propositions and theorems stated in this section are deferred to Appendix \ref{apdx:exactVI} and \ref{apdx:R-LSVI}.\\
In the inexact VI defined in Procedure \ref{prd:inexact_VI}, the updates take the form
\begin{equation}
    \hat{P}_{i+1} =\mathcal{R} (\hat{P}_{i}, S) + \Delta_{i}
    =\mathcal{R} (\hat{P}_{i}, S + \Delta_{i})
\end{equation}
where $\mathcal{R}$ is the Riccati operator defined by
\begin{equation}
\label{eq:riccati-operator}
    \mathcal{R}\left(P,S\right) \triangleq A^{\top} P A - A^{\top} P B(R+B^{\top} P B)^{-1} B^{\top} P A+ S 
\end{equation}
for $P\in\Sset^{n}_{+}$.
The operator $\mathcal{R}$ plays a central role in our analysis. Its structure allows us to reduce the study of $\hat{P}_{i+1}$’s dependence on $\Delta_i$ to understanding how the Riccati solution depends on the cost matrix $S$.

\subsection{Convergence of the Exact VI}\label{sec:convergence-of-the-exact-VI}
\noindent We begin by establishing the convergence of exact VI in this section, laying the groundwork for proving the convergence and robustness of inexact VI.

The following theorem shows that exact VI is globally exponentially stable, thereby extending earlier results, which either establish local exponential stability \cite{song2024convergence} or guarantee exponential stability only when the initial value function lies within the semidefinite cone of the optimal solution \cite{Lee2022TAC,zhang2023revisiting}. 

\begin{thm}[Global Exponential Convergence of Exact VI]
\label{thm:KLexactVI}
The exact value iteration in Procedure \ref{prd:exact_VI} is globally exponentially stable. That is, there exists a $\mathcal{KL}$-function $\beta$, such that $\forall i \in \Zset_+ \text{ and }  \forall P_{0} \in \Sset_+^{n}$
\begin{equation}
\label{eq:expGlobal}
    \norm{P_i - P^*}_{2} \le \beta(\norm{P_0 - P^*}_{2}, i).
\end{equation}
In addition, $\beta(r, i) = {\alpha}\theta^{2i}r$ for some constant ${\alpha} \ge 1$ and $\theta\in(\rho(A^{*}),1)$.
\end{thm}

The next proposition asserts that a sequence of value matrices generated by iterating on the DARE will dominate another such sequence (in the $\succeq $ sense), if the cost matrices associated with the first sequence dominate those of the second. This result is the discrete analogue of applying Grönwall’s inequality \cite{teschl2024ordinary} to the matrix differential equation $\dot{P}(t) = \mathcal{R}(P(t), S(t))$.

\begin{prop}\label{prop:comparisionRDE}
Let $\{P_i^j\}_{i\in \Zset_+}$ with $j\in\{1,2\}$ be the following sequences
\begin{equation}
    P_{i+1}^{j} = \mathcal{R} (P_{i}^{j}, S_{i}^{j} ), \quad P_{0}^{j} \in\Sset^{n}_{+}.
\end{equation}
Suppose that $S^{1}_i \succeq S^{2}_i$ for all $i \in \Zset_+$ and $P_{0}^{1} \succeq P_{0}^{2}$. Then, $P_{i}^{1} \succeq P_{i}^{2}$ for all $i \in \Zset_+$.
\end{prop}

Proposition \ref{prop:continuouslyDiffP} states that the unique solution to the Riccati equation with perturbed cost matrix  $S + \delta I_n$  depends smoothly on the perturbation parameter $\delta$.

\begin{prop}\label{prop:continuouslyDiffP}
    Let $\delta > - \eigmin{S}$ and $P^*(\delta)\in\Sset^{n}_{++}$ be the unique solution of the following DARE:
    \begin{equation}
        \label{eq:riccati_delta}
        P =\mathcal{R}\left(P,  S + \delta I_{n}\right).
    \end{equation}
    Then $P^*(\delta)$ is continuously differentiable for all $\delta \in (- \eigmin{S}, \infty)$.
\end{prop}

\subsection{Convergence of the Inexact VI}\label{sec:convergence-of-the-inexact-VI}
\noindent This section focuses on the analysis of inexact VI in Procedure \ref{prd:inexact_VI}, which serves as a general framework for examining the convergence properties of VI for a broad range of disturbances. Using Proposition \ref{prop:comparisionRDE}, we first show that inexact VI can be bounded between two instances of exact VI. Then, leveraging Proposition \ref{prop:continuouslyDiffP}, we demonstrate that these bounds can be made arbitrarily tight in a continuous manner, concluding that inexact VI is small-disturbance ISS.

\begin{thm}[Small-Disturbance ISS of Inexact VI]
\label{thm:inexact_VI}
    Inexact value iteration is small-disturbance ISS. Specifically, if $\norm{\Delta}_{\infty} < \eigmin{S}$, there exist a $\mathcal{KL}$-function $\beta$, and a $\mathcal{K}$-function $\gamma$, such that 
    \begin{equation}
        \norm{\hat{P}_i - P^*}_{2} \le \beta(\norm{\hat{P}_0 - P^*}_{2}, i) + \gamma(\norm{\Delta}_{\infty}) 
    \end{equation}
    for all $i \in \Zset_+$ and all $\hat{P}_0 \in\Sset^{n}_{+}$.
\end{thm}

\begin{cor}
    \label{cor:compactC}
    For a given $\hat{P}_{0}\in\Sset_{+}^{n}$, let $\mathcal{C}(\hat{P}_{0})$ denote the set of value matrices generated by the inexact value iteration in Procedure \ref{prd:inexact_VI}, 
    starting from $\hat{P}_{0}$, with all possible $\left\{\Delta_{i}\right\}_{i=0}^{\infty}$ satisfying $\|\Delta \|_{\infty}<\eigmin{S}$. Then the closure $\bar{\mathcal{C}}(\hat{P}_{0})$ is compact.
\end{cor}

\begin{rem}
Theorem \ref{thm:inexact_VI} provides a global result that extends the local exponential asymptotic stability results of continuous-time systems shown in \cite{bian2019continuous} and the exponential asymptotic stability for any $\hat{P}_{0}\succeq P^{*}$ established in \cite{Lee2022TAC,zhang2023revisiting}.
\end{rem}

\subsection{Convergence of R-LSVI}
\noindent \label{sec:convergence-of-the-robust-VI}
In this section, we establish the convergence of the R-LSVI algorithm by demonstrating that it is a special case of inexact VI. Disturbances in the R-LSVI algorithm (Algorithm \ref{alg:R-LSVI}) result from approximation errors in $Q_{i} = \mathcal{Q}(\hat{P}_{i})$ of the form $\hat{Q}_{i} =  Q_{i}+\Delta Q_{i}$, which are due to insufficient data, system noise, and numerical inaccuracies. 
Theorem \ref{thm:robust_LSVI} establishes the convergence of R-LSVI.

\begin{prop}
\label{prop:robust_VI}
    For any $\hat{P}_0 \in \mathbb{S}_+^n$, there exists $d^{*}>0$, a $\mathcal{KL}$-function $\beta$ and a $\mathcal{K}$-function $\tilde{\gamma}$, such that if $\norm{\Delta Q}_{\infty} < d^{*}$, then
    \begin{equation}
        \norm{\hat{P}_i - P^*}_{2} \le \beta(\norm{\hat{P}_0 - P^*}_{2}, i) + \tilde{\gamma}(\norm{\Delta Q}_{\infty})
    \end{equation}
    for all $i \in \Zset_+$.
\end{prop}

The following proposition shows that there exists $T_{0}>0$, independent of the iteration step $i$, such that for any sample size $T>T_{0}$, the output of Algorithm \ref{alg:R-LSVI} satisfies $\norm{\Delta Q_{i}}_{\infty}<d^{*}$.

\begin{prop}
    \label{prop:robust_LSVI}
    For a given $\hat{P}_{0}$, let $\mathcal{P}(\hat{P}_{0})$ denote the set of value matrices generated by Algorithm \ref{alg:R-LSVI}, starting from $\hat{P}_{0}$, with all possible $\left\{\Delta Q_{i}\right\}_{i=0}^{\infty}$ satisfying $\|\Delta Q\|_{\infty}<d^{*}$ with $d^{*}$ defined in Proposition \ref{prop:robust_VI}. Then, under Assumption \ref{assum:invertible1} there exist $T_{0}\in\Zset^{+}$ such that, for any $T>T_{0}$,
    \begin{equation}
        \hat{P}_{i}\in\mathcal{P}(\hat{P}_{0})\implies \|\Delta Q_{i} \|_{2}<d^{*}, \text{ almost surely}.
    \end{equation}
\end{prop}
We proceed to prove the convergence of R-LSVI.
\begin{thm}[Convergence of R-LSVI]
\label{thm:robust_LSVI}
Under Assumption \ref{assum:invertible1}, for Algorithm \ref{alg:R-LSVI}, given any initial value matrix $\hat{P}_{0}\in\Sset^{n}_{+}$ and any $\epsilon>0$, there exists $T_{0} \in \Zset_{+}$ such that for all $T \geq T_{0}$ 
$$
\limsup _{i \rightarrow \infty} \norm{ \hat{P}_{i}-P^{*}}_{2}<\epsilon \text{ almost surely} \,.
$$
\end{thm} 
\begin{proof}
Since $\hat{P}_{0} \in \mathcal{P}(\hat{P}_{0})$, by Proposition \ref{prop:robust_LSVI}, $\left\|\Delta Q_{1}\right\|_{2}<d^{*}$ almost surely. By definition, $\hat{P}_{1} \in \mathcal{P}(\hat{P}_{0})$. Thus, by induction, we have $\left\|\Delta Q\right\|_{\infty}<d^{*}$ almost surely. By Proposition \ref{prop:robust_VI}, the result follows.
\end{proof}

\section{Empirical Results}\label{sec:empirical-results}
\noindent We examine the convergence, stability, and robustness of the R-LSVI algorithm using a standard RL benchmark  problem with the unstable Laplacian dynamics, given by \eqref{eq:dynamic}--\eqref{eq:functional} with 
\begin{equation} 
\label{data:1}
    A = \left[\begin{array}{lll} 1.01 & 0.01 & 0\\ 0.01 & 1.01 & 0.01\\ 0 & 0.01 & 1.01\\ \end{array}\right], 
    \quad B = \left[\begin{array}{lll} 1 & 0& 0\\ 0 & 1 & 0\\ 0 & 0 & 1\\ \end{array}\right]
\end{equation}
\begin{equation}
\label{data:2}
    C = S = I_3 , \quad R = 1000I_3\,.
\end{equation}
Several authors have adopted this problem, which represents a stylized ``data center cooling'' scenario \cite{gao2014machine}, to evaluate different RL approaches for stochastic LQR \cite{tu2018least,abbasi2019model,dean2020sample}.

\begin{rem}
The open-loop system given by the matrices in \eqref{data:1} and \eqref{data:2} is unstable. Specifically, for any nonzero initial condition, the open-loop state sequence diverges. Should an RL algorithm (explicitly or implicitly) underestimate a diagonal entry of $A$ as less than one, it will erroneously treat the corresponding mode as stable and allocate insufficient control effort, especially given the high control cost represented by $R = 1000I_3$. Consequently, obtaining a high-quality estimate of the system's true behavior is imperative for achieving near-optimal control.
\end{rem}

% For comparison, we include the following well-established methods: Nominal Control \cite{recht2019tour,dean2020sample}, Policy Gradient \cite{montague1999reinforcement}, LSPI \cite{krauth2019finite}, and O-LSPI \cite{pang2021robust}. Details on these procedures can be found in the Appendix~\ref{apdx:procedure_comparison}.
For comparison, we include the following well-established methods:
\begin{itemize}
    \item \textbf{Nominal Control} \cite{recht2019tour,dean2020sample}: a plug-in PI (VI) model-based procedure that first estimates the system and cost matrices $(A,B,S,R)$ and then solves the stochastic LQR problem using Hewer’s PI \cite{hewer1971iterative} (Bertsekas's VI; see \cite{bertsekas2011dynamic}).
    \item \textbf{Policy Gradient} \cite{montague1999reinforcement}: a model-free approach that directly learns a policy function from episodic experiences. In our implementation, we use the Adam optimizer \cite{kingma2014adam} for policy gradient updates and estimate the cost using the running cost over ten consecutive points. We supply the algorithm with the exact gradient. 
    \item \textbf{LSPI} \cite{krauth2019finite}: a classical PI-based ADP method. For a fair comparison with R-LSVI, we adapt LSPI to incorporate the additive noise term $\mu(P)$ and to include the cost matrices $(S,R)$ explicitly in the Hamiltonian $\mathcal{Q}(P)$.
    \item \textbf{O-LSPI} \cite{pang2021robust}: an ADP procedure consisting of two loops, an outer PI loop and an inner VI loop. We set $I_{\text{max}}^{\text{inner}} = 20$ and $I_{\text{max}}^{\text{outer}} = 5$ to enhance the overall accuracy of the policy updates.
\end{itemize}
For fairness of comparison, all methods use the same reset-based sampling $\eqref{eq:dyn_reset1}, \eqref{eq:dyn_reset2}$ and the same rescaling by $\alpha_{t}^{2}$.
\begin{rem}
    Each method yields a control gain $\hat{K}$. If $\hat{K}$ is stabilizing, the linear state feedback control policy $u_{t}=-\hat{K}x_{t}$ is applied to the system, resulting in the cost $J\left(u\right)=\Tr{C^{\top} \hat{P} C}$, where the value matrix $\hat{P}\in\Sset^{n}_{++}$ is the unique  positive definite solution of the Lyapunov equation
\begin{equation}
    \label{eq:lyapunov}
    \hat{P} = (A-B\hat{K})^{\top} \hat{P}(A-B\hat{K})+S+\hat{K}^{\top}R\hat{K}\,.
\end{equation}
This distinction matters for value iteration because the controller is derived from $\hat{P}$ (computed via the Lyapunov equation \eqref{eq:lyapunov}) rather than from $\hat{P}_{I_{\text{max}}}$.
\end{rem}
For each problem instance and method, we run the algorithms $100$ times and calculate the fraction of trials that produce a stabilizing control gain $\hat{K}$. To ensure comparability, all algorithms are executed for a total of $100$ iterations $\left(I_{\text{max}} = 100\right)$ and use the same number of data points.

We initialize each procedure with $x_{\text{init}} = 0_{3}$ and $\hat{P}_{0} = \beta I_{3}$, with $\beta\sim\mathcal{U}\left(0,1\right)$ to ensure that $\hat{P}_{0}\in\Sset^{n}_{+}$. Since the nominal PI, LSPI and O-LSPI are PI-based algorithms, it is necessary to select an initial control gain, $\hat{K}_0$. For consistency and fairness across comparisons, we choose $\hat{K}_0$ as the corresponding induced gains $\hat{K}_0=(R+B^{\top}\hat{P}_0B)^{-1}B^{\top}\hat{P}_0A$, as in \cite{song2024convergence}. For the policy gradient, we choose a  stabilizing  initialization $\hat{K}_{0}  = \alpha_{0} I_{3}$, with $\alpha_{0}\sim\mathcal{U}\left(0.1,1\right)$ to enhance its convergence.

\subsection{Convergence Analysis}
\noindent We evaluate algorithmic convergence on the benchmark problem as the sample size increases. In our experiments, we assess convergence using the cost relative error, defined as 
\begin{equation}
  \frac{\Tr{C^{\top}(\hat{P}-P^{*})C}}{\Tr{C^{\top}P^{*}C} } 
\end{equation}
where the optimal $P^{*}$ is given by formula \eqref{eq:riccati-eqn}. In addition, we compute the fraction of instances in which the controller $\hat{K}$ is stable. These measures allow us to precisely quantify  the performance of the algorithms relative to the optimal solution. To address the non-stabilizing data-collection policy, we set ${K}_c = \alpha I_{3}$, with $\alpha\sim\mathcal{U}\left(-0.1,0\right)$, take $ \Sigma_{\eta} = I_{m}$ and fix the reset bound to $d=1000$. Finally, to examine the effect of the rescaling step, we evaluate R-LSVI both with and without the rescaling factor $\alpha_{t}$.
\begin{figure}[!ht]
    \centering
    \hspace*{-0.10cm}
    \includegraphics[scale=\globalscale]{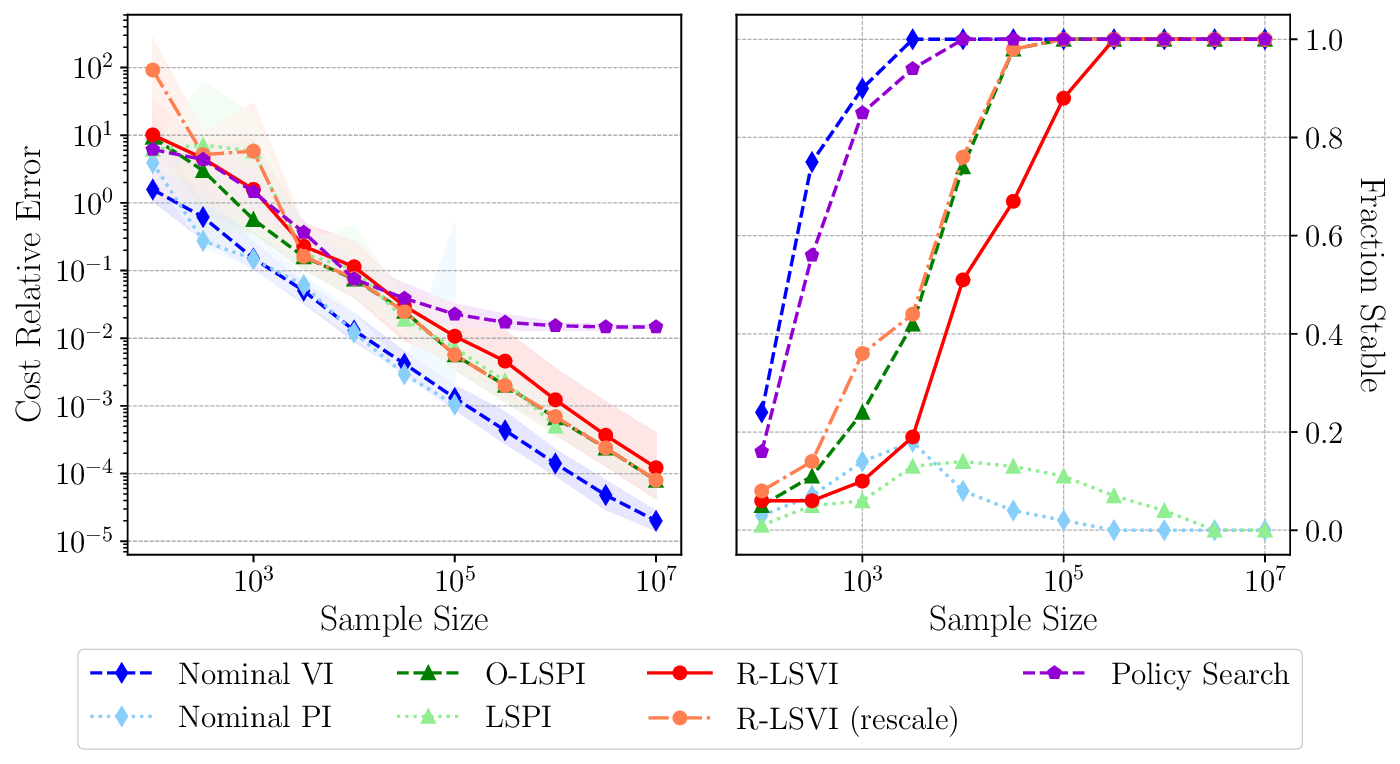}
    \caption{Convergence and stability of the nominal control (PI and VI), O-LSPI, LSPI, R-LSVI (with and without rescaling) and policy gradient algorithms on the data center cooling benchmark problem as the sample size increases. (Left) Cost relative error. Dashed lines represent the median relative error, with the shaded region covering the $25^{th}$ to $75^{th}$ percentiles, estimated from $100$ trajectories. (Right) Frequency of stabilizing controllers found by the algorithms. Only costs corresponding to stable control gains are plotted in the left panel.}
    \label{fig:convergence}
\end{figure}
Figure~\ref{fig:convergence} illustrates the convergence behavior as a function of the sample size. The left panel  displays the convergence of the cost relative error, while the right panel shows the fraction of stable controllers. Our numerical experiments show that all VI algorithms eventually converge to the optimal solution, with a cost relative error on the order of $10^{-5}$, whereas the policy gradient method attains a cost relative error on the order of $10^{-2}$. However, only VI-based procedures achieve full stabilization at sample size of about $10^{5}$, whereas PI-based procedures do not produce stabilizing controllers. This highlights that VI does not require a stabilizing initial policy. Despite being a PI method, O-LSPI behaves like VI due to its inner VI loop.\\
We observe that the rescaling step plays a central role in stabilizing the method under uneven data collection. Indeed, R-LSVI attains full stability for a sample size of $10^{4}$, whereas its non-rescaling version requires about $10^{5}$. Note that nominal VI converges slightly faster than the ADP methods, consistent with \cite{krauth2019finite,abbasi2019model,tu2018least,tu2019gap}.\\
Finally, the policy gradient baseline benefits from a favorable initialization $\hat{K}_{0}$ and access to exact gradients; therefore, it serves as an optimistic benchmark for model-free methods. Nevertheless, R-LSVI matches or exceeds its performance.

For the remainder of our empirical analysis, we adopt a stabilizing data-collection policy that allows us to focus on other aspects of our procedure. Under this policy, the sampled data are more homogeneous and concentrate around zero, eliminating the need for the reset and rescaling steps.

\subsection{Adaptivity under Non-Quadratic Costs}
\noindent We present an example illustrating the adaptability of the algorithms to changes in the cost function for the benchmark problem. Specifically, we use the same system dynamics \eqref{data:1}--\eqref{data:2}, but change the cost function to
\begin{equation}
\label{eq:non-quad cost}
\tilde{c}(x_{t}, u_{t}) \triangleq x_{t}^{\top}Sx_{t} + (|u_{t}|^{\frac{\kappa}{2}})^{\top}R|u_{t}|^{\frac{\kappa}{2}} 
\end{equation}
where $|\cdot|^{\frac{\kappa}{2}}$ is applied elementwise and $\kappa \in [1,3]$. Since this cost function is non-quadratic in  $u$ for values other than $\kappa=2$, the optimal cost is not expected to adhere to the standard form $\Tr{C^{\top}P^{*}C}$. Consequently, direct comparisons of relative costs are not meaningful in this setting. Instead, we evaluate the cost obtained using a control input $u_t = -\hat{K} x_t$, where $\hat{K}$ is computed for each method using a sample size of $T = 10^5$. The choice of $T$ is guided by prior results (Figure~\ref{fig:convergence}), which show that this sample size enables most methods to reach $100\%$ stability. To assess adaptability to cost misspecification independently of the non-stabilizing data-collection policy, we set $K_c=\alpha I_{3}$ with $\alpha\sim\mathcal{U}(0.1,0.2)$.
\begin{figure}[!ht]
    \centering
    \hspace*{-0.10cm}
    \includegraphics[scale=\globalscale]{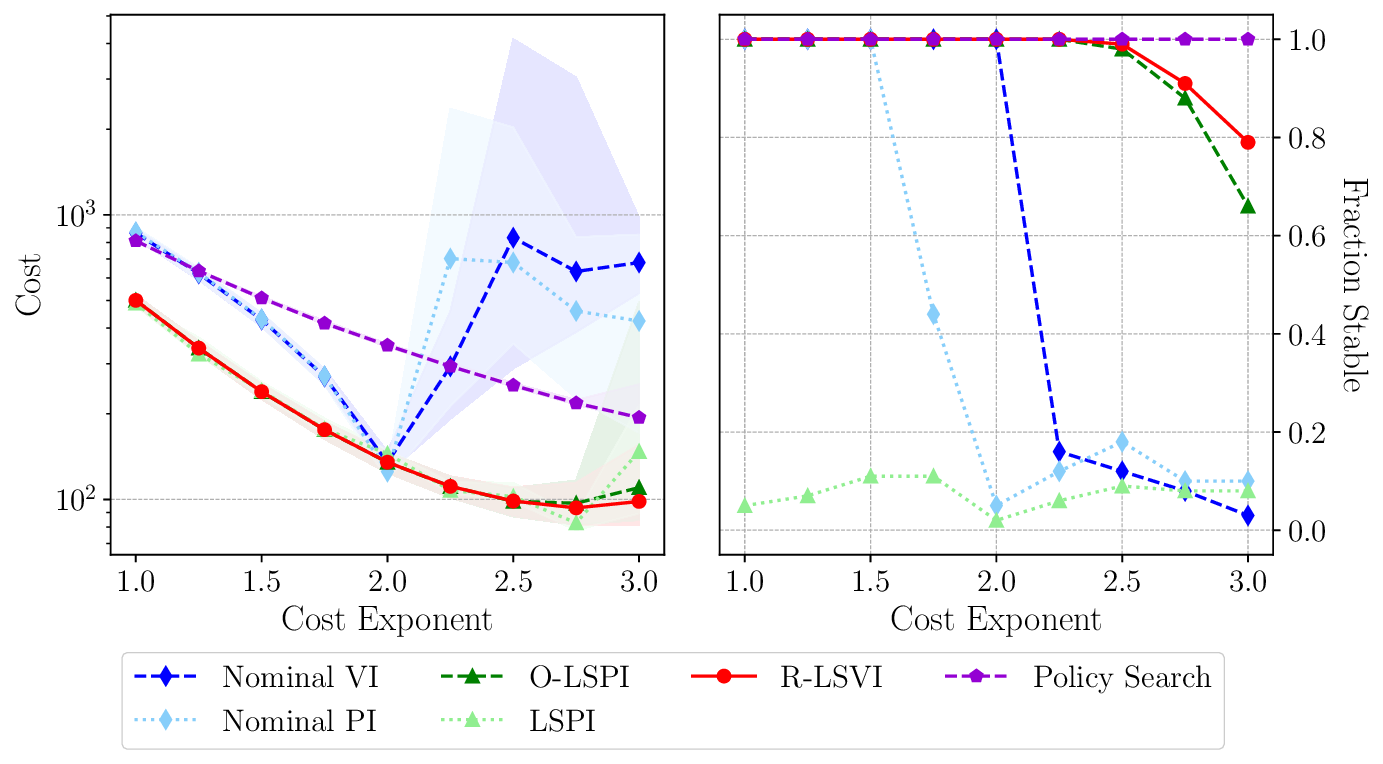}
    \caption{Convergence and stability of the nominal control (VI and PI), 
  LSPI, O-LSPI, R-LSVI and policy gradient algorithms on the data center cooling benchmark problem with non-quadratic cost. (Left) Costs as the exponent $\kappa$ varies. Dashed lines represent the median relative error, with the shaded region covering the $25^{th}$ to $75^{th}$ percentiles, estimated from $100$ trajectories. (Right) Frequency of stabilizing controllers found by the algorithms. Only costs corresponding to stable control gains are plotted in the left panel.}
    \label{fig:adaptability}
\end{figure}
In Figure~\ref{fig:adaptability}, we observe that the ADP algorithms consistently outperform nominal model-based control in both adaptability and stability. In the left panel, for $\kappa=2$, the nominal control slightly outperforms the other methods, which is consistent with the exact model specification and the findings in Figure~\ref{fig:convergence}. However, beyond this point, R-LSVI systematically outperforms all other algorithms, exhibiting superior adaptability. Notably, the right panel shows that for $\kappa > 2$, nominal model-based control becomes unstable, with stability dropping sharply below  $20\%$. In contrast, non-model-based methods maintain robust performance, with R-LSVI achieving the highest stability, remaining above $80\%$ across all values of $\kappa$. 
\begin{rem}
\label{rem:policy-grad}
We emphasize that, theoretically, the policy gradient method enjoys a significant advantage because it has direct access to the modified cost function and the ability to compute the exact gradient through backpropagation. Thus, it has the best chance of identifying the implicit optimal linear control compared to the other methods. Nonetheless, we observe that in this example R-LSVI delivers the best controllers in terms of cost.
\end{rem}

\section{An Application to Dynamic Portfolio Allocation}
\label{sec:portfolio-allocation}
\noindent In this section, we apply the R-LSVI algorithm to a practical financial problem, referred to as dynamic portfolio allocation, using a modeling framework similar to G\^{a}rleanu and Pedersen \cite{garleanu2013dynamic}. This model admits a standard stochastic LQR representation \cite{abeille2016lqg, moallemi2017dynamic}. The state vector $x_{t}\in\mathbb{R}^{2N + M}$ is the concatenation of the portfolio weights at time $t-1$ ($w_{t-1} \in \mathbb{R}^N$), mean-reverting strategy signals ($f_t \in \mathbb{R}^M$) and asset returns between $[t-1,t)$ ($r_{t} \in \mathbb{R}^N$), while the control input $u_{t}\in\Rset^{N}$  represents the change in portfolio weights at time $t$ ($u_{t} \triangleq w_{t} - w_{t-1}$). In this model we assume the asset returns, $r_{t+1}$, are a linear combination of the mean-reverting signal $f_t$ leading to the linear dynamics 
\begin{equation}
\label{eq:matrix_gp1}
    A \triangleq \left[\begin{array}{ccc} I_{N} & 0 & 0 \\ 0 & I_{M}-\Phi & 0\\ 0 & \Pi & 0 \end{array}\right], 
    B \triangleq \left[\begin{array}{c} I_{N} \\ 0 \\0 \end{array}\right], 
    C \triangleq \left[\begin{array}{cc} 0 & 0\\ \Omega^{\frac{1}{2}}&0\\0&\Sigma^{\frac{1}{2}} \end{array}\right]
\end{equation}
where $\Pi\in\Rset^{N\times M}$ is the matrix of factor loadings, $\Phi\in\Rset^{M\times M}$ is the matrix of mean-reversion coefficients, with $\rho(I_{M}-\Phi)<1$, and $\Sigma\in\Sset^{N}_{++}$ and $\Omega\in\Sset^{M}_{++}$ are the covariance matrices of the white noise processes $\{\epsilon^{r}_{t}\}$ and $\{\epsilon^{f}_{t}\}$, representing return noise and factor innovations, respectively.\\
The investor aims to select a dynamic trading strategy $(u_{t})_{t>0}$ that maximizes long-term average return adjusted for risk (defined as the portfolio variance) after price impact, leading to the following cost matrices
\begin{equation}
\label{eq:matrix_gp2}
    S \triangleq \left[\begin{array}{ccc} \gamma \Sigma & 0 & -I_{N} \\ 0 & \star_{f} & 0\\ -I_{N} & 0 & \star_{r} \end{array}\right], R \triangleq \Lambda
\end{equation}
where $\gamma>0$ is the risk aversion coefficient and $\Lambda\in\Sset^{M}_{++}$ is the matrix that quantifies price impact.

We remark that adding quadratic terms in $f_{t}$ and $r_{t}$ to the cost function does not affect the optimal control since these variables are not controllable. This gives us flexibility in choosing the matrices $\star_{f}\in\Rset^{M\times M}$ and $\star_{r}\in\Rset^{N\times N}$ so that $\eigmin{S}>0$ is satisfied. In particular, as we demonstrate in the Appendix \ref{apdx:GP}, setting $\star_{f} = I_{M}$ and $\star_{r} = 2\gamma^{-1}\Sigma^{-1}$ meets this requirement. 

In our model, the factors $\{f_{t}\}$ serve as predictors for forecasting future returns.  In real-world trading applications, developing such factors, referred to as ``alphas,'' is a crucial element of successful dynamic portfolio allocation. While the development of predictors is beyond the scope of this article, we simulate real-world alphas by intentionally introducing look-ahead bias and defining one factor for each asset as
\begin{equation}
    f_{i,t} \triangleq \frac{1}{100}\sum_{k=0}^{99} r_{i, t+1-k}\,.
\end{equation}
This choice of factors yields linear dynamics with a predictive $R^2$ of about $1\%$, which is generally considered a meaningful alpha signal in practice. In the quantitative finance literature, this approach of simulating alphas is known as bootstrap or synthetic alphas \cite{kolm2022mean}.

In this article, we set $N=3$ and choose the assets to be the three commodity futures contracts for coffee, cocoa and sugar from January 1, 1996 to January 3, 2024. We estimate the covariance matrices $\Sigma$ and $\Omega$
using returns over the full sample, shrinking the correlations $50\%$ toward zero as in \cite{garleanu2013dynamic}. We take $\gamma = 30$, $\Lambda = 0.03 I_{3}$, $x_{0} = 0_{9}$, $K_{c} = 0_{3 \times 9}$ and $\Sigma_{\eta} = I_{3}$.
\begin{figure}[!ht]
    \centering
    \hspace*{-0.10cm}
    \includegraphics[scale=\globalscale]{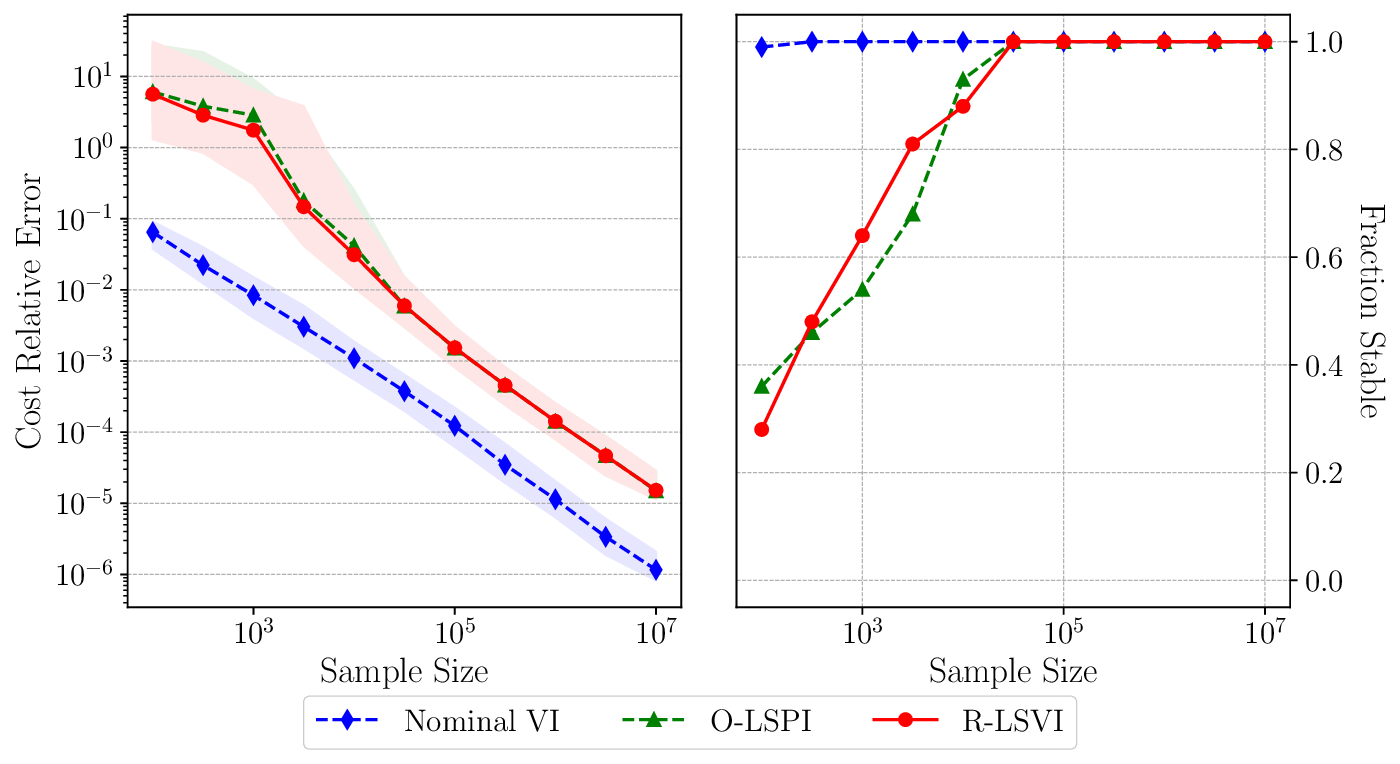}
    \caption{Convergence and stability of the nominal VI, O-LSPI and R-LSVI algorithms on the dynamic portfolio allocation problem as the sample size increases. (Left) Cost relative error. Dashed lines represent the median relative error, with the shaded region covering the $25^{th}$ to $75^{th}$ percentiles, estimated from $100$ trajectories. (Right) Frequency of stabilizing controllers found by the algorithms. Only costs corresponding to stable control gains are plotted in the left panel.}
    \label{fig:convergence_gp}
\end{figure}
As in Section \ref{sec:empirical-results}, we conduct two experiments for the dynamic portfolio allocation problem: (a) convergence as a function of sample size, and (b) adaptability to non-quadratic cost misspecification. 
We exclude the policy gradient method and compare R-LSVI, O-LSPI, and the model-based nominal control. Overall, the resulting convergence and stability behavior is qualitatively similar to that observed in the data center cooling example.
Specifically, the left panel of Figure~\ref{fig:convergence_gp} demonstrates that all methods converge as the sample size increases, with the non-model-based procedure achieving a cost relative error near $10^{-5}$, and nominal control exhibiting slightly faster convergence. Meanwhile, the right panel shows that all methods attain $100\%$ stability when the sample size reaches approximately $10^4$.

\begin{figure}[!ht]
    \centering
    \hspace*{-0.10cm}
    \includegraphics[scale=\globalscale]{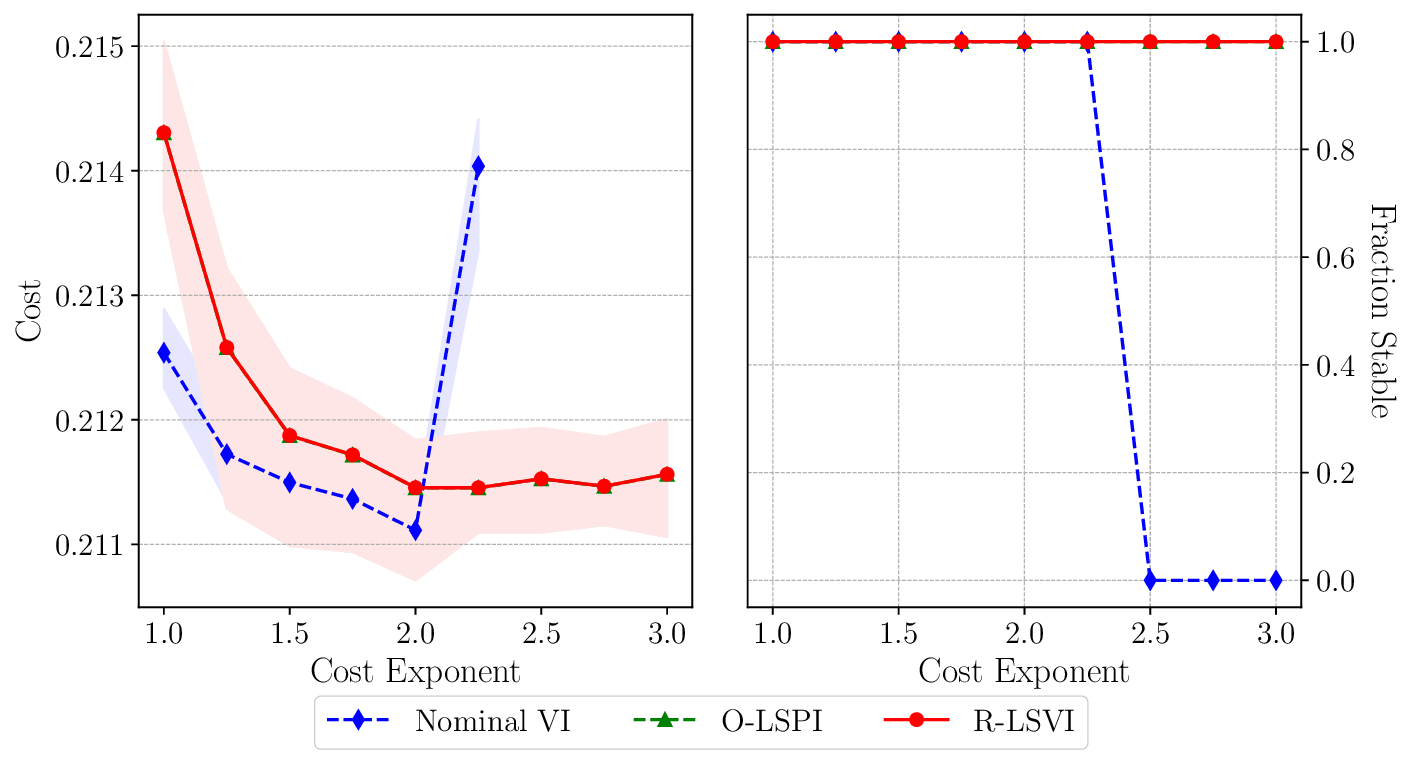}
    \caption{Convergence and stability of the nominal VI, O-LSPI and R-LSVI algorithms on the dynamic portfolio allocation problem with non-quadratic cost. (Left) Costs as the exponent $\kappa$ varies. Dashed lines represent the median relative error, with the shaded region covering the $25^{th}$ to $75^{th}$ percentiles, estimated from $100$ trajectories. (Right) Frequency of stabilizing controllers found by the algorithms. Only costs corresponding to stable control gains are plotted in the left panel.}
    \label{fig:adaptability_gp}
\end{figure}

Regarding cost misspecification, the left panel of Figure~\ref{fig:adaptability_gp} shows that the nominal model-based control slightly outperforms the non-model-based methods for $\kappa \leq 2$,  in contrast to the previous example, where this occurred only at $\kappa = 2$. In the right panel, however, when $\kappa > 2$, the nominal approach fails to produce a stable control gain, leading to a dramatic increase in final cost. In contrast, the non-model-based methods maintain robust performance, achieving $100\%$ stability with final costs of comparable magnitude.

\section{Conclusions}\label{sec:conclusions}
\noindent In this article, we established the convergence, robustness, and stability of value iteration (VI) for data-driven control of stochastic linear quadratic (LQ) systems in discrete time with completely unknown dynamics and cost matrices. 
%Existing ADP methods for discrete-time systems with continuous state spaces primarily rely on policy iteration (PI), which requires an initial stabilizing policy. In contrast, we develop the first off-policy, robust, non-model-based VI algorithm, termed R-LSVI, which eliminates the need for a prior stabilizing control policy.
By leveraging small-disturbance input-to-state stability (ISS), we prove that exact VI is globally exponentially stable and that inexact VI is robust to external disturbances, in the sense that the generated sequence remains within a small neighborhood of the optimal solution when the disturbances are sufficiently small.

To validate our theoretical findings, we evaluated R-LSVI against established non-model-based, model-based, and model-free methods using the ``data center cooling'' benchmark problem. Our empirical results revealed that R-LSVI consistently outperformed other non-model-based approaches in both stability and convergence speed. 
In addition, R-LSVI demonstrated strong adaptability to non-quadratic costs, underscoring its potential for broader control applications.  Finally, we presented a fully data-driven solution to dynamic portfolio allocation, an important problem in quantitative finance, where our non-model-based approach consistently demonstrated reliable convergence and strong adaptability to non-quadratic costs.

\section*{Acknowledgments}
\noindent Gr\'egoire G.~Macqueron gratefully acknowledges support from the Jeffrey and Diane Rosenbluth Fellowship.

\appendix

% \section{Procedure Comparison}\label{apdx:procedure_comparison}
% In Section~\ref{sec:empirical-results} of the main article, we include the following well-established methods:
% \begin{itemize}
%     \item \textbf{Nominal Control} \cite{recht2019tour,dean2020sample}: a plug-in PI (VI) model-based procedure that first estimates the system and cost matrices $(A,B,S,R)$ and then solves the stochastic LQR problem using Hewer’s PI \cite{hewer1971iterative} (Bertsekas's VI; see \cite{bertsekas2011dynamic}).
%     \item \textbf{Policy Gradient} \cite{montague1999reinforcement}: a model-free approach that directly learns a policy function from episodic experiences. In our implementation, we use the Adam optimizer \cite{kingma2014adam} for policy gradient updates and estimate the cost using the running cost over ten consecutive points. We supply the algorithm with the exact gradient. 
%     \item \textbf{LSPI} \cite{krauth2019finite}: a classical PI-based ADP method. For a fair comparison with R-LSVI, we adapt LSPI to incorporate the additive noise term $\mu(P)$ and to include the cost matrices $(S,R)$ explicitly in the Hamiltonian $\mathcal{Q}(P)$.
%     \item \textbf{O-LSPI} \cite{pang2021robust}: an ADP procedure consisting of two loops, an outer PI loop and an inner VI loop. We set $I_{\text{max}}^{\text{inner}} = 20$ and $I_{\text{max}}^{\text{outer}} = 5$ to enhance the overall accuracy of the policy updates.
% \end{itemize}

\section{Auxiliary Results}\label{apdx:aux_results}
Throughout this supplement, we use the following notation. For all $i\in\{\emptyset,\Zset_{+}\}$ and $j \in \{*,1,2\}$, define
\begin{equation}
        K_i^j \triangleq (R+B^\top P_{i}^{j} B)^{-1}B^\top P_{i}^{j}A, \quad
        A_{i}^j \triangleq A - BK_{i}^j, \quad
        R_{i}^j \triangleq R+B^\top P_{i}^{j} B\,.
\end{equation}

\begin{lem}
\label{lem:kron} \
    \begin{enumerate}[(i)]
    \item For any $A \in \Rset^{n\times m}, B \in \Rset^{m\times p}$ and $C \in \Rset^{p\times q}$, we have 
    \begin{equation}
    \label{prop:kron1} 
    \vect{ABC}=\left(C^{\top} \otimes A\right) \vect{B} \in\Rset^{np} \,.
    \end{equation}
    \item For any $v_1 \in \Rset^{n}$ and $v_2 \in \Rset^{m}$, we have 
    \begin{equation}
    \label{prop:kron2}
    v_{1} \otimes v_{2}=\vect{v_{2}v_{1}^{\top}} \in\Rset^{nm} \,.
    \end{equation}
    \item For any $v_1 \in \Rset^{n}, v_2 \in \Rset^{m}$ and $M \in \Rset^{n\times m}$, we have 
    \begin{equation}
    \label{prop:kron3}
    v_{1}^{\top} M v_{2}=\left(v_2^{\top} \otimes v_1^{\top}\right) \vect{M} \in\Rset \,.
    \end{equation}
    \item For $v \in \Rset^{n}$ and $S \in  \Sset^{n}$, we have 
    \begin{align}
    v^{\top} S v
    = \vect{vv^{\top}}^{\top} \vect{S}
    = \tilde{v}^{\top} \svec{S} \,.
    \end{align}
\end{enumerate}
\end{lem}
\begin{proof}
    See \cite{magnus2019matrix}.
\end{proof}

\begin{lem}
\label{lem:lin_alg_trick}
For any $X \in \Rset^{m \times n}$ and $Y \in \Rset^{n \times m}$ such that both $(I_m + XY)$ and $(I_n + YX)$ are invertible, we have
\begin{enumerate}[(i)]
    \item $(I_m + XY)^{-1}X = X(I_n + YX)^{-1}\,.$
    \item $(I_m + XY)^{-1} = I_{m} - (I_m + XY)^{-1}XY\,.$
    \item $I_m - X (I_n+YX)^{-1}Y = (I_m+ XY)^{-1}\,.$
\end{enumerate}
\end{lem}
\begin{proof}
\begin{align}
    {(i)} &= (I_m + XY)^{-1}X (I_n + YX)(I_n + YX)^{-1}
    \nonumber \\&
    = X(I_n + YX)^{-1} \,.\\
    {(ii)} &= I_{m} - (I_m + XY)^{-1}(I_m + XY) + (I_m + XY)^{-1}
     \notag \\& 
    = I_{m} - (I_m + XY)^{-1}XY \,.\\
    {(iii)} &= I_m - (I_m+XY)^{-1}XY
    \notag \\&
    = (I_m + XY)^{-1} \,.
\end{align}
\end{proof}

\begin{lem}
\label{lem:Schur_power_bound}
    Let $A\in\Rset^{n\times n}$ such that $\rho(A)<1$, then for any $\theta\in (\rho(A),1)$ there exists a constant $C_{\theta}\geq 1$ such that
    \begin{equation}
        \norm{A^{i}}_{2}\leq C_{\theta}\theta^{i},\quad \forall i\geq 0
    \end{equation}
\end{lem}
\begin{proof}
    See \cite[Lemma~5.6.10]{horn2012matrix}
\end{proof}

\begin{lem}
    \label{lem:Delta_Q_positive}
    For any $P\in\Sset^{n}_{+}$, we denote by $\hat{\mathcal{Q}}\left(P\right) \triangleq \mathcal{Q}\left(P\right) + \Delta \mathcal{Q}\left(P\right)$, where $\Delta \mathcal{Q}\left(P\right) \in \Sset^{n+m}$ is a disturbance, and define $d \triangleq \min(\eigmin{S}, \eigmin{R})$. The following properties hold
    \begin{enumerate}[(i)]
        \item $ \mathcal{Q}\left(P\right)\succeq  \mathcal{Q}\left(0_{n\times n}\right)\succ  0$.
        \item If $\norm{\Delta \mathcal{Q}\left(P\right)}_{2} < d$ then $\hat{\mathcal{Q}}\left(P\right)\in\Sset^{n+m}_{++}$.
        \item If $\norm{\Delta \mathcal{Q}\left(P\right)}_{2} < d$ then $[\hat{\mathcal{Q}}\left(P\right)]_{uu}$ is invertible.
    \end{enumerate}
\end{lem}
\begin{proof}
{(i)} We have that $\mathcal{Q}\left(0_{n\times n}\right) = \begin{bmatrix} S & 0 \\ 0 & R \end{bmatrix}\succ  0$,
as $S\in\Sset^{n}_{++}$ and $R\in\Sset^{m}_{++}$. Hence, as $P\in\Sset^{n}_{+}$ we have
\begin{equation}
    \mathcal{G}\left(P\right) = \mathcal{Q}\left(P\right) -   \mathcal{Q}\left(0_{n\times n}\right) = \begin{bmatrix} A \\ B \end{bmatrix}^{\top}P\begin{bmatrix} A \\ B\end{bmatrix}\succeq 0\,.
\end{equation}
{(ii)} Using {(i)}, we have that for any $P\in\Sset^{n}_{+}$
\begin{align}
   \eigmin{\hat{\mathcal{Q}}(P)} 
   \geq  d - \norm{\Delta \mathcal{Q}\left(P\right)}_{2} > 0.
\end{align}

{(iii)} Recall that $[\hat{\mathcal{Q}}\left(P\right)]_{uu} = R + B^{\top}PB + [\Delta\mathcal{Q}\left(P\right)]_{uu}$, hence
\begin{equation}\eigmin{[\hat{\mathcal{Q}}\left(P\right)]_{uu}}
\geq  \eigmin{R} -\norm{\Delta \mathcal{Q}\left(P\right)}_{2}
> 0
\end{equation}
as $\eigmin{R}>0$. 
\end{proof}
\begin{lem}\label{lem:Rpositive}
    For any $P\in\Sset^{n}_{+}$, $R\in\Sset^{m}_{++}$ and $B\in\Rset^{n\times m}$, we have
    \begin{equation}
    \label{eq:P-PBBPPSD}
       P \succeq P - PB (R+B^\top P B)^{-1}B^\top P \succeq 0 \,.
    \end{equation} 
\end{lem}
\begin{proof} The first inequality follows directly from $PB (R+B^\top P B)^{-1}B^\top P \in\Sset^{n}_{+}$. For the second inequality, we denote by $M =P^{\frac{1}{2}} B$ and $N = MR^{-\frac{1}{2}}$. Then, 
\begin{equation}
\label{eq:PBBPPSD}
        M (R+M^\top M)^{-1}M^\top = N(I_m +N^\top N)^{-1}N^\top
        \preceq I_n\,.
    \end{equation}
By multiplying inequality \eqref{eq:PBBPPSD} by $P^{\frac{1}{2}}$ from both left and right, the proof is complete.
\end{proof}

\begin{lem}\label{lem:Kdiff}
    Suppose that $P^i, P^j \in\Sset^{n}_{+}$.
    Then, 
    \begin{equation}
        K^i - K^j = R^{-1}B^\top \Lambda_i^{-\top}(P^i - P^j) A^j 
    \end{equation}
    where $\Lambda_i = I_n + BR^{-1}B^\top P^i$.
\end{lem}

\begin{proof}
By the definition of $K^j$, 
\begin{align}
    K^j &= R^{-1}(I_m + B^\top P^j BR^{-1})^{-1}B^\top P^j A
    = R^{-1}B^\top P^j \Lambda_j^{-1}A.
\end{align}
Hence, we have
\begin{align}
    K^i -  K^j &= R^{-1}B^\top (P^i - P^j) \Lambda_j^{-1}A 
    \\&\quad\notag
    +  R^{-1}B^\top P^i \Lambda_i^{-1}A- R^{-1}B^\top P^i \Lambda_j^{-1}A \,.
\end{align}
Let us focus on the second term $R^{-1}B^\top P^i \Lambda_i^{-1}A- R^{-1}B^\top P^i \Lambda_j^{-1}A\triangleq(*)$
\begin{align}
    (*) & = R^{-1}B^\top P^i (\Lambda_i^{-1} \Lambda_j \Lambda_j^{-1} - \Lambda_i^{-1} \Lambda_i\Lambda_j^{-1})A \\
    &=R^{-1}B^\top P^i \Lambda_i^{-1}( \Lambda_j -  \Lambda_i)\Lambda_j^{-1}A \\
    &= - R^{-1}B^\top P^i \Lambda_i^{-1} BR^{-1}B^\top (P^i - P^j) \Lambda_j^{-1}A \,.
\end{align}
Finally,
\begin{align}
    K^i -  K^j &= R^{-1}B^\top[ I_n - P^i \Lambda_i^{-1} BR^{-1}B^\top] (P^i - P^j) \Lambda_j^{-1}A 
    \notag\\&
    = R^{-1}B^\top \Lambda_i^{-\top} (P^i - P^j) \Lambda_j^{-1}A
\end{align}
where the last equality follows from Lemma \ref{lem:lin_alg_trick}~{(iii)}. We conclude by observing that $A^{j} = {\Lambda_j}^{-1}A$ from Lemma \ref{lem:lin_alg_trick}~{(iii)}.
\end{proof}

\begin{lem}\label{lem:Pdiff}
Let $\{P_i^j\}_{i\in \Zset_+}$ with $j\in\{1,2\}$ be the sequences of the following Riccati difference equations
\begin{equation}
\label{eq:twoRDE}
    P_{i+1}^{j} = \mathcal{R} (P_{i}^{j}, S_{i}^{j} ), \quad P_{0}^{j} \in\Sset^{n}_{+}.
\end{equation}
Then the sequence $\Delta{P}_i = P_{i}^1 - P_{i}^2$ follows the following dynamic
\begin{align}\label{eq:PdiffRDE1}
    \Delta{P}_{i+1} &= (A_i^2)^\top [\Delta{P}_{i}
    - \Delta{P}_{i} B (R+B^\top {P}^{1}_i B)^{-1}B^\top \Delta{P}_{i} ]A_i^2 + \Delta S_{i}
\end{align}
with $\Delta S_{i} = S^{1}_{i} - S^{2}_{i}$.
\end{lem}

\begin{proof}
    Equation \eqref{eq:twoRDE} with $j=1$ can be rewritten as
\begin{align}\label{eq:P1RDE}
        P_{i+1}^1 &= A^\top P_{i}^1 A + S^{1}_i - (K_{i}^1)^\top R_{i}^1K_i^1\\
        &= (A_i^2)^\top P_{i}^1 A_i^2 +S^{1}_i - (K_{i}^1)^TR_{i}^1K_i^1
        + (K_{i}^2)^\top B^\top P_{i}^{1}A + A^{\top}P_{i}^{1}BK_{i}^2
        \notag\\&\qquad 
        - (K_i^2)^\top (B^\top P_{i}^{1} B) K_{i}^2\\
        &= (A_i^2)^\top P_{i}^1 A_i^2 +S^{1}_i - (K_{i}^1)^TR_{i}^1K_i^1
        + (K_{i}^2)^\top R_{i}^{1}K_{i}^{1}+ (K_{i}^1)^\top R_{i}^{1}K_{i}^2
        \notag\\&\qquad 
        - (K_i^2)^\top (R_{i}^{1}-R) K_{i}^2 \,.
\end{align}
    In addition, we can rewrite \eqref{eq:twoRDE} with $j=2$ as 
    \begin{equation}
    \label{eq:P2RDE}
        P_{i+1}^2 = (A_i^2)^\top P_{i}^2 A_i^2 + S^{2}_i +(K_i^2)^\top R K_{i}^2 \,.
    \end{equation}
    Subtracting \eqref{eq:P2RDE} from \eqref{eq:P1RDE} yields
    \begin{align}\label{eq:PdiffRDE2}
        \Delta{P}_{i+1} &= (A_i^2)^\top \Delta{P}_i A_i^2 + \Delta S_{i}
        - (K_i^1 - K_i^2)^\top (R+B^\top P_{i}^{1} B)(K_i^1 - K_i^2).
    \end{align}
    We conclude by applying Lemma \ref{lem:Kdiff} to \eqref{eq:PdiffRDE2}.
\end{proof}

\begin{lem}[Matrix Differentiation, \cite{magnus2019matrix}]
\label{lem:mat_diff}
Let $M\in\mathbb{R}^{n \times n}$ and $S \subset \mathbb{R}^{n \times q}$ be an open subset such that for every $P\in S$, $(I_{n} + PM)$ is invertible. Define the matrix-valued function $U: S \rightarrow \mathbb{R}^{n \times n}$ by $U(P) \triangleq (I_{n} + PM)^{-1}P$.
Then $U$ is differentiable on $S$ and the Jacobian is the $n^{2} \times n^{2}$ matrix
\begin{equation}
    \frac{d\vect{U(P)}}{d\vect{P}} =(I_n + MP)^{-\top}\otimes (I_n + PM)^{-1} \,.
\end{equation}
\end{lem}

\begin{proof}
We note that for any matrix function $F: \mathbb{R}^{n \times p} \rightarrow \mathbb{R}^{m \times q}$, $d\vect{F(X)} = \vect{dF(X)}$. Then
\begin{align}
d\left((I_n+PM)^{-1}P\right) 
\notag&
    =  d\left((I_n+PM)^{-1}\right) P + (I_n+PM)^{-1}d\left(P\right)\\
    \notag
    &=  (I_n+PM)^{-1}d\left(P\right)  (I_n+PM)^{-1}
    + (I_n+PM)^{-1}d\left(P\right)\\
    \notag
    &=  (I_n+PM)^{-1}d\left(P\right)\left[I_n - M(I_n+PM)^{-1}P\right]\\
    &= (I_n+PM)^{-1}d\left(P\right)(I_n+MP)^{-1}
\end{align}
where the last two equalities follow from Lemma \ref{lem:lin_alg_trick}~{(iii)}. By applying $\vect{\cdot}$ on both sides of the identity above, the proof is complete.
\end{proof}

\section{Exact and Inexact VI}\label{apdx:exactVI}
\begin{lem}\label{lem:Ppositive}
    If $\norm{\Delta}_\infty < \eigmin{S}$, then $\hat{P}_i \succeq 0$ for all $i \in \Zset_+$. 
\end{lem}

\begin{proof}
The proof is by induction. Note that $\norm{\Delta}_\infty = \sup\{\norm{
    \Delta_i}
    , i\in \Zset_+\} < \eigmin{S}$ implies that $S+\Delta_i \succ 0$ for all $i \in \Zset_+$. Suppose $\hat{P}_i \succeq 0$ for some $i \in \Zset_+$. Since
    \begin{align}
        \hat{P}_{i+1}&=\mathcal{R}(\hat{P}_{i}, S + \Delta_i)
        =A^{\top}\left(\hat{P}_{i} - \hat{P}_{i}B (R+B^\top \hat{P}_{i} B)^{-1}B^\top \hat{P}_{i}\right)A + S + \Delta_i
    \end{align}
    it follows that $\hat{P}_{i+1} \succeq 0$ by Lemma \ref{lem:Rpositive}. As $\hat{P}_0 \succeq 0$, the proof is complete.
\end{proof}

\begin{lem}\label{lem:localExp}
Consider the exact value iteration
\begin{equation}
    {P}_{i+1} = \mathcal{R}\left({P}_{i}, S\right),\quad {P}_0 \in\Sset^{n}_{+}.   
\end{equation}
Fix any $\theta\in(\rho(A^*),1)$ and let $C_{\theta}$ be the constant from Lemma~\ref{lem:Schur_power_bound}. Set 
\begin{equation}
    \nu \triangleq \norm{A^*}_2^{2}\norm{BR^{-1}B^\top}_2,
    \qquad
    \delta_\theta \triangleq \frac{\theta^{2}(1-\theta^{2})}{4\nu C_{\theta}^{4}}.
\end{equation}
If $\Delta P_0 \triangleq P_0 - P^* \in \mathcal{P} \triangleq \{P \in \mathbb{S}^{n}: \norm{P}_{2} \le \delta_\theta\}$, then for all $i\in \Zset_+$
\begin{equation}
\label{eq:strong-induction}
    \norm{\Delta P_i}_{2} \le 2C_{\theta}^{2}\theta^{2i} \norm{\Delta P_0}_{2}.
\end{equation}
\end{lem}
\begin{proof}
We prove \eqref{eq:strong-induction} by strong induction. For $i=0$, the claim is immediate since $2C_{\theta}^{2}\geq 1$.\\
Fix $i>0$ and suppose \eqref{eq:strong-induction} holds for all $j\leq i$. Then observing that $\Delta P_{i} = P_{i} - P^{*}$, applying Lemma~\ref{lem:Pdiff} and Lemma~\ref{lem:lin_alg_trick}~{(iii)}, we have
\begin{align}
    \Delta{P}_{i+1} 
    &= {(A^*)^\top \Delta{P}_{i}}[ I_{n} - B (R+ B^{\top}{P}_iB)^{-1} B^{\top}\Delta{P}_{i} ]A^*\\
    \label{eq:PPoptDiff0}
    &= (A^*)^\top \Delta{P}_{i} 
    \times [ I_{n} -  BR^{-1}B^{\top} (I_{n}+{P}_i BR^{-1}B^{\top})^{-1} \Delta{P}_{i} ]A^*.
\end{align}
Decompose~\eqref{eq:PPoptDiff0} into
\begin{equation}
\label{eq:delta_pi_epsilon_i}
    \Delta P_{i+1} = (A^*)^\top \Delta P_i A^* + \epsilon_i,
\end{equation}
where
\begin{equation}
\label{eq:epsilon_i_def}
    \epsilon_i
    \triangleq
    -(A^*)^\top \Delta P_i
    \big[BR^{-1}B^\top (I_n + P_i BR^{-1}B^\top)^{-1}\big]
    \Delta P_i A^*.
\end{equation}
Hence,
\begin{equation}
\label{eq:epsilon_i_bound}
    \norm{\epsilon_i}_{2}
    \le
    \norm{A^*}_{2}^{2}\norm{BR^{-1}B^\top}_{2}\norm{\Delta P_i}_{2}^{2}
    = \nu \norm{\Delta P_i}_{2}^{2}.
\end{equation}
Unrolling \eqref{eq:delta_pi_epsilon_i} yields
\begin{equation}\label{eq:unroll}
    \Delta P_{i+1} = \big((A^*)^{i+1}\big)^\top \Delta P_0 (A^*)^{i+1} 
    + \sum_{j=0}^{i}\big((A^*)^{i-j}\big)^\top \epsilon_j (A^*)^{i-j}.
\end{equation}
Taking the norm of~\eqref{eq:unroll} and applying Lemma~\ref{lem:Schur_power_bound} to each power of $A^*$ gives
\begin{equation}\label{eq:keybound}
    \norm{\Delta P_{i+1}}_2 \le C_{\theta}^{2} \theta^{2(i+1)} \norm{\Delta P_0}_2 
    + \nu C_{\theta}^{2}\sum_{j=0}^{i}\theta^{2(i-j)} \norm{\Delta P_j}_2^{2}.
\end{equation}
By strong induction we have that $\norm{\Delta P_j}_2 \le 2C_{\theta}^{2}\theta^{2j}\norm{\Delta P_0}_2$ for all $j\leq i$, hence
\begin{align}
\notag
    \norm{\Delta P_{i+1}}_2 
    &\le C_{\theta}^{2}\theta^{2(i+1)}\norm{\Delta P_0}_2 + 4\nu C_{\theta}^{6}\norm{\Delta P_0}_2^{2}\,\theta^{2i}\sum_{j=0}^{i}\theta^{2j}.
\end{align}
Using $\sum_{j=0}^{i}\theta^{2j} \le 1/(1-\theta^{2})$ we obtain
\begin{equation}
    \norm{\Delta P_{i+1}}_2 \le C_{\theta}^{2}\theta^{2(i+1)}\norm{\Delta P_0}_2\!\left[1 + \frac{4\nu C_{\theta}^{4}\norm{\Delta P_0}_2}{\theta^{2}(1-\theta^{2})}\right].
\end{equation}
Since %$\Delta P_0 \in \mathcal{P}$, i.e. 
$\norm{\Delta P_0}_2 \le \delta_\theta$, the bracket is bounded by $2$, which closes the induction.
\end{proof}

\subsection{Proof of Proposition \ref{prop:comparisionRDE}}
    The proof is by induction. Clearly, $P_{i}^{1} \succeq P_{i}^{2}$ for $i=0$. Assume $P_{i}^{1} \succeq P_{i}^{2}$ for some $i \in \Zset_+$. Using Lemma \ref{lem:Pdiff} and that $R + B^{\top}P^{1}_{i}B = R^2_i + B^\top \Delta{P}_{i} B$, we obtain 
\begin{align}
        \Delta{P}_{i+1} &= \Delta S_{i} + (A_i^2)^\top  [ \Delta{P}_i
        -  \Delta{P}_i B(R^2_i + B^\top \Delta{P}_{i} B)^{-1} B^\top \Delta{P}_i]A_i^2.
    \end{align}
By applying Lemma \ref{lem:Rpositive}, we have that $\Delta{P}_{i+1} \succeq 0$, i.e.~$P_{i+1}^1 \succeq P_{i+1}^2$.  

\subsection{Proof of Proposition \ref{prop:continuouslyDiffP}}
The proof follows from the implicit function theorem. We observe that  from Lemma \ref{lem:lin_alg_trick} we have
\begin{align}
    \mathcal{R}\left(P, S\right) 
    \notag &= A^\top [I_{n} -  P BR^{-1}(I_{m}+B^\top PBR)^{-1}B^\top ] P A + S\\
    &= A^\top(I_n+PBR^{-1}B^\top)^{-1}PA + S
\end{align}
Noting that we can rewrite equation \eqref{eq:riccati_delta} as $ 0 =\mathcal{R}\left(P,  S + \delta I_{n} - P\right)$, define 
\begin{align}\label{eq:partialF}
    f(\vect{P} ; \delta) & \triangleq \vect{\mathcal{R}\left(P, S + \delta I_{n} - P\right)} \\
    &= A^\top \otimes A^\top \vect{(I_n+PBR^{-1}B^\top)^{-1}P}
    \notag\\  &\qquad
    + \vect{S+\delta I_n} -  \vect{P} \,.
\end{align} 
Let $M = BR^{-1}B^\top$. Using Lemma \ref{lem:mat_diff} and that $(I_n + MP)^{-\top} = (I_n + PM)^{-1}$ for $P$, $M \in \Sset^{n}$, we have
\begin{equation}
    \frac{d \vect{(I_n+PM)^{-1}P}}{d\vect{P}} =(I_n + PM)^{-1}\otimes (I_n + PM)^{-1} \,.
\end{equation}
Hence,
\begin{equation}
\label{eq:partialF2}
    \frac{d f(\vect{P}; \delta)}{d \vect{P}} = A(P)^\top \otimes A(P)^\top - I_{n^2}
\end{equation}
where $A(P) \triangleq (I_n+ MP)^{-1}A$. We apply Lemma \ref{lem:lin_alg_trick}~{(iii)} to rewrite $A(P)$ as 
\begin{equation}
    A(P) = A - BK(P)
\end{equation}
where $K(P) = (R+B^TPB)^{-1}B^\top PA$. Therefore, for $P = P^*(\delta)$ with \linebreak$\delta \in (-\eigmin{S}, \infty)$, $A\left(P^*(\delta)\right)$ is Schur stable (see, \cite[Section 2.4]{lewis2012optimal} for details). Since $\rho(A\left(P^*(\delta)\right)) < 1$, the left hand side of \eqref{eq:partialF2} is invertible. Consequently, by the implicit function theorem $P^*(\delta)$ is continuously differentiable with respect to $\delta$.

\subsection{Proof of Theorem \ref{thm:KLexactVI}}
Let $\Delta P_i = P_i - P^*$. Then, using Lemma \ref{lem:Pdiff} and that $R + B^{\top}P_{i}B = R^{*} + B^\top \Delta{P}_{i} B$, we obtain
\begin{align}
    \label{eq:DeltaPpos}
    \Delta{P}_{i+1} &= (A^*)^\top  M_{i}A^* %\,.
\end{align}
with 
\begin{equation}
    M_i\triangleq \Delta P_i-\Delta P_iB(R^*+B^\top\Delta P_iB)^{-1}B^\top\Delta P_i\,.
\end{equation}
The statement in the theorem is obviously true if $\norm{A^*}_{2} = 0$. Next, we assume $\norm{A^*}_{2} \neq 0$ and consider the statement in three cases: \textcircled{1} $\Delta{P}_{0}$ is positive semi-definite, i.e. $P_{0} \succeq P^{*}\succeq 0$; \textcircled{2} $\Delta{P}_{0}$ is negative semi-definite, i.e. $0\preceq P_{0} \preceq P^{*}$; and \textcircled{3} $\Delta{P}_{0}$ is indefinite, i.e. $\eigmax{\Delta{P}_{0}} > 0$ and $\eigmin{\Delta{P}_{0}} < 0$.
\\
\textcircled{1} $\Delta P_{0} \succeq 0$: First, we show by induction that $\Delta P_{i}\succeq 0$ for all $i\in\Zset_{+}$. This is true at $i=0$ by assumption. Suppose that $\Delta P_{i}\succeq 0$ for some $i\in\Zset_{+}$. Then Lemma~\ref{lem:Rpositive} yields $0\preceq M_i\preceq\Delta P_i$, and thus
$0\preceq\Delta P_{i+1}=(A^*)^\top M_iA^*\preceq (A^*)^\top\Delta P_iA^*$. Hence $\Delta P_{i+1}\succeq 0$. Iterating the upper bound gives
$0\preceq\Delta P_i\preceq ((A^*)^i)^\top\Delta P_0(A^*)^i$. Applying Lemma~\ref{lem:Schur_power_bound}, we obtain
\begin{equation}
    \norm{\Delta P_i}_{2} \le C_{\theta}^{2}\theta^{2i}\norm{\Delta P_0}_{2}\,.
\end{equation}
\\
\textcircled{2} $\Delta{P}_{0} \preceq 0$: 
 Let $\{ P_i^0\}_{i \in \Zset_+}$ denote the sequence generated by the exact VI with the initial condition $P_0^0 = 0$, and $\Delta P_i^0 = P_i^0 - P^*$. By \cite[Proposition 4.4.1]{bertsekas2011dynamic}, $\Delta P_i^0$ is nondecreasing with respect to $i$ and $\lim_{i \to \infty} \Delta P_i^0 = 0$. Let $i_{0}$, independent of $P_{0}$, denote the smallest index $i$ such that $\Delta P_i^0 \in \mathcal{P}$. Since $\Delta P_0^0 \preceq \Delta P_0 \preceq 0$, by  Proposition~\ref{prop:comparisionRDE}, $\Delta P_{i_{0}}^0 \preceq \Delta P_{i_{0}} \preceq 0$.
%\\
 Hence, $\norm{\Delta P_{i_{0}}}_{2} \le \norm{\Delta P_{i_{0}}^0}_{2}$, and so $\Delta P_{i_{0}} \in \mathcal{P}$. Therefore, by Lemma \ref{lem:localExp}, for all $i \ge i_{0}$, we have 
\begin{equation}
\label{eq:normDeltaP1}
    \norm{\Delta P_i}_{2} \le \frac{2C_{\theta}^{2}}{\theta^{2i_{0}}}\theta^{2i} \norm{\Delta P_0}_{2}\,.
\end{equation}
For any $0 \le i \le i_{0}$, we have
\begin{align}\label{eq:normDeltaP2}
    \notag
    \norm{\Delta P_i}_{2} &\le \norm{P^*}_{2} \\ 
    \notag
    &\le \norm{P^*}_{2} \frac{ \norm{\Delta P_0}_{2}}{\delta_{\theta}}
    \\ 
    &\le \frac{\norm{P^*}_{2} }{\delta_{\theta}\theta^{2i_{0}}} \theta^{2i} \norm{\Delta P_0}_{2} %\,.
\end{align}
with $\delta_{\theta}$ from Lemma~\ref{lem:localExp}.
It follows from \eqref{eq:normDeltaP1} and \eqref{eq:normDeltaP2} that 
\begin{equation}
    \norm{\Delta P_i}_{2} \le {\gamma}\theta^{2i} \norm{\Delta P_0}_{2}, \quad \forall i \in \Zset_{+} %.
\end{equation}
with $\gamma = \max\left\{\frac{2C_{\theta}^{2}}{\theta^{2i_{0}}}, \frac{\norm{P^*}_{2} }{\delta_{\theta}\theta^{2i_{0}}}\right\}$.
\\
\textcircled{3} $\Delta{P}_{0}$ is indefinite: Define
\begin{equation}
\Delta P^{*}_{0} \triangleq (P^*)^{-\frac{1}{2}}\Delta P_{0}(P^*)^{-\frac{1}{2}}
\end{equation}
and $\eta \triangleq -\eigmin{\Delta P^{*}_{0}}$. Since $P_{0}\succeq 0$, we have $\Delta P_{0}\succeq -P^{*}$, hence by congruence $\Delta P^{*}_{0}\succeq -I_{n}$, which yields $\eta \in [0,1]$.
\\
Let $\{P_{i}^{+}\}_{i\in\Zset_{+}}$ and $\{P_{i}^{-}\}_{i\in\Zset_{+}}$ be the sequences generated by the exact VI from $P_{0}^{+}=\norm{\Delta P_{0}}_{2}I_{n}+P^{*}$ and $P_{0}^{-}=(1-\eta)P^{*}$, respectively. Since $\eta\in[0,1]$, both $P_{0}^{+}$ and $P_{0}^{-}$ lie in $\Sset^{n}_{+}$.
\\
Then we have $P_{0}^{-}\preceq P_{0}\preceq P_{0}^{+}$. The upper inequality is immediate. The lower inequality is equivalent to $\Delta P_{0}+\eta P^{*}\succeq 0$, which holds by congruence with $\Delta P^{*}_{0}-\eigmin{\Delta P^{*}_{0}} I_{n}\succeq 0$.
\\
Since $\Delta P_{0}^{+}=\norm{\Delta P_{0}}_{2}I_{n}\succeq 0$, Case \textcircled{1} gives
\begin{equation}
    \norm{\Delta P_{i}^{+}}_{2}\le  C_{\theta}^{2}\theta^{2i}\norm{\Delta P_{0}}_{2}.
\end{equation}
Since $\Delta P_{0}^{-}=-\eta P^{*}\preceq 0$, Case \textcircled{2} gives
\begin{equation}
    \norm{\Delta P_{i}^{-}}_{2}\le \gamma\theta^{2i}\eta\norm{P^{*}}_{2}.
\end{equation}
Furthermore, we have $\eta \le \frac{\norm{\Delta P_{0}}_{2}}{\eigmin{P^{*}}}$, hence
\begin{equation}
    \norm{\Delta P_{i}^{-}}_{2}\le \gamma\,\kappa(P^{*})\,\theta^{2i}\norm{\Delta P_{0}}_{2}
\end{equation}
with $\kappa(P^{*})\triangleq \frac{\norm{P^{*}}_{2}}{\eigmin{P^{*}}}$.
\\
Finally, by Proposition~\ref{prop:comparisionRDE} we have for all $i$, $\Delta P_{i}^{-}\preceq \Delta P_{i}\preceq \Delta P_{i}^{+}$ and $\Delta P_{i}^{-}\preceq 0\preceq \Delta P_{i}^{+}$. Hence $\norm{\Delta P_{i}}_{2}\le \max\{\norm{\Delta P_{i}^{-}}_{2},\norm{\Delta P_{i}^{+}}_{2}\}$, so
\begin{equation}
    \norm{\Delta P_{i}}_{2}\le \max\{C_{\theta}^{2},\gamma\kappa(P^{*})\}\theta^{2i}\norm{\Delta P_{0}}_{2}\,.
\end{equation}
To conclude, combining the three cases, we obtain
\begin{equation}
\norm{\Delta P_{i}}_{2}\le \alpha\theta^{2i}\norm{\Delta P_{0}}_{2} = \beta(\norm{P_0 - P^*}_{2},i)
\end{equation}
 with $\alpha = \gamma\kappa(P^{*})$ and the proof is complete.

\subsection{Proof of Theorem \ref{thm:inexact_VI}}\label{apdx:inexactVI}
    Define $\delta^+ \triangleq \norm{\Delta}_{\infty}$, $\delta^- \triangleq -\norm{\Delta}_{\infty}$, and consider the sequences $\{ P_i^+\}_{i\in \Zset_{+}}$ and $\{ P_i^-\}_{i\in \Zset_{+}}$ generated by
\begin{align}
   {P}_{i+1}^+ &= \mathcal{R}\left({P}_{i}^+, S + \delta^+I_n\right), \quad {P}^+_0 = P_0 \\
   {P}_{i+1}^- &= \mathcal{R}\left({P}_{i}^-, S + \delta^-I_n\right), \quad {P}^-_0 = P_0.
\end{align}    
Since the pair $(A,B)$ is stabilizable and $S + \delta^+I_n \succ 0$ and $S + \delta^-I_n \succ 0$, both ${P}^+_{i}$ and ${P}^-_{i}$ converge to ${P}_+^{*}$ and ${P}_-^{*}$, which are the unique positive-definite solutions of
\begin{align}
    {P}_{+}^* &= \mathcal{R}\left({P}_{+}^*, S + \delta^+I_n\right)\\
    {P}_{-}^* &= \mathcal{R}\left({P}_{-}^*, S + \delta^-I_n\right).
\end{align}    
By Proposition \ref{prop:continuouslyDiffP}, there exists an $L > 0$ such that $\norm{P_+^* - P^*}_{2} \le L\norm{\Delta}_{\infty}$ and $\norm{P_-^* - P^*}_{2} \le L\norm{\Delta}_{\infty}$. By Theorem \ref{thm:KLexactVI} and the weak triangle inequality \cite{jiang1994small} we obtain
\begin{align}
    \norm{P_i^+ - P_+^*}_{2} &\le \beta(\norm{P_0 - P_+^*}_{2},i) 
    \le \beta(\norm{P_0 - P^*}_{2},i) +  \gamma(\norm{\Delta}_{\infty})\\
    \norm{P_i^- - P_-^*}_{2} &\le \beta(\norm{P_0 - P_-^*}_{2},i)
    \le \beta(\norm{P_0 - P^*}_{2},i)+  \gamma(\norm{\Delta}_{\infty}).
\end{align}    
By Proposition \ref{prop:comparisionRDE} we have 
\begin{equation}
    P_i^- - P^* \preceq \hat{P}_i - P^* \preceq P_i^+ - P^*.
\end{equation}
Hence, we have
\begin{align}
    \norm{\hat{P}_i - P^*}_{2} &\le \max\{ \norm{P_i^+ - P^*}_{2}, \norm{P_i^- - P^*}_{2}\} 
    \notag\\&
    \le \beta(\norm{P_0 - P^*}_{2},i)+  \gamma(\norm{\Delta}_{\infty})
\end{align}
and the proof is complete.

\section{R-LSVI}\label{apdx:R-LSVI}

\subsection{Proof of Proposition \ref{prop:reset}}
\label{apdx:reset}
Controllability: Consider the non-resetting chain
\begin{equation}
    X_{t+1} = LX_t + D\xi_{t+1}, \quad \xi_t \underset{iid}{\sim}\mathcal N(0_{m+p}, I_{m+p})
\end{equation}
with $L = A - BK_c$ and $ D = (B\Sigma_\eta^{1/2} ~ C)$. Without loss of generality, we assume that $(L,D)$ is controllable. Since $X_0=x_0 = 0_n$, all uncontrollable state components 
are identically zero for all $t$. Their invariant distribution is therefore the point mass $\delta_0$.\\
Irreducibility: Define the $t$-step reachability Gramian
\begin{equation}
    W_t = \sum_{s=0}^{t-1} L^s D D^\top (L^s)^\top.
\end{equation}
Since $(L,D)$ is controllable, there exists $T\le n$ such that $\mathrm{rank}(W_T)=n$ and hence $W_T\succ0$.
Let $\bar{\mathcal B}_1(0_n)$ be the closed unit ball and $\mathcal V\subset \bar{\mathcal B}^{\infty}_d(0_n)$ any non-empty open set, define
\begin{equation}
    E(\mathcal V) \triangleq \{X_0 = 0_{n},\ X_T \in \mathcal V,\ X_i \in \bar{\mathcal B}^{\infty}_d(0_n)\ \forall i<T\}.
\end{equation}
Because the transition density of $(X_t)$ is Gaussian and $W_T\succ0$, we have 
$\mathbb P(E(\mathcal V))>0$.
Since $(x_{t})_{t>0}$ and $(X_{t})_{t>0}$ coincide while inside $\bar{\mathcal B}^{\infty}_d(0_n)$, the resetting chain $(x_{t})_{t>0}$ can reach any non-empty open set of $\bar{\mathcal B}^{\infty}_d(0_n)$ in a finite number of steps $T\le n$ from $0_{n}$ with strictly positive probability.\\
Positive Recurrence: Let $p \triangleq \mathbb P(\|D\xi_t\|_2 > 2d)$.
Since $DD^\top\ne 0$, we have $p>0$, independent of $t$ and the current state $x_{t}$.
Let $\tau = \inf\{t\ge0:\|X_t\|_{\infty} > d\}$.
At each step, the probability of exiting the ball (and thus resetting to $0$) is at least $p$, hence
\begin{equation}
    \mathbb P(\tau > k) \le (1-p)^k \text{ and } \mathbb E[\tau] \le \frac{1}{p} < \infty.
\end{equation}
Thus the reset mechanism defines a regenerative Markov chain with finite mean cycle length and atom $\{0_n\}$, implying positive recurrence.\\
Hence, the regenerative Markov chain admits a unique invariant probability measure $\pi$.
\begin{equation}
       \pi(A)\;=\; \frac{\mathbb{E}\left[\sum_{t=0}^{\tau}\mathbbm{1}_{\{x_t\in A\}}\right]}{\mathbb{E}[\tau]}, \quad A\subset \bar{\mathcal{B}}^{\infty}_{d}(0_{n})\,.
\end{equation}

\subsection{Proof of Proposition \ref{prop:robust_VI}}\label{apdx:proof_prop_vi}
    Consider the set $\mathcal{C}(\hat{P}_{0})$ of $\Sset^{n}_{+}$ defined in Corollary \ref{cor:compactC}. By Lemma \ref{lem:Delta_Q_positive}, we know that for $\norm{\Delta Q_{i}}_{2} < d$, $[Q_{i}]_{uu}$ and $[\hat{Q}_{i}]_{uu}$ are invertible. Therefore, we relate $\Delta_{i}$ and $\Delta Q_{i}$ via 
    \begin{equation}
        \Delta_{i} =  \mathcal{H}(\hat{Q}_{i}) - \mathcal{H}\left(Q_{i}\right)
    \end{equation}
    where $\mathcal{H}$ is an analytic function. The set $\{Q_{i} = \mathcal{Q}(\hat{P}_{i}) | \hat{P}_{i}\in\bar{\mathcal{C}}(\hat{P}_{0})\}$ is compact as it is the image of a compact set under the continuous function $\mathcal{Q}$. By the extreme value theorem there exists a constant $c(\bar{\mathcal{C}}(\hat{P}_{0}))>0$ such that 
    \begin{equation}
        \norm{\Delta_{i}}_{\infty} \leq c(\bar{\mathcal{C}}(\hat{P}_{0}))\norm{\Delta Q_{i}}_{\infty} \,.
    \end{equation}
    Therefore, with $d^{*} \triangleq \min \left(d, \frac{\eigmin{S}}{2c(\bar{\mathcal{C}}(\hat{P}_{0}))}\right)$, the sequence $\{\Delta_{i}\}_{i\geq0}$ generated by Algorithm \ref{alg:R-LSVI} satisfies 
    \begin{equation}
        \left\|\Delta_{i}\right\|_{\infty}<\eigmin{S}
    \end{equation}
    for all $i>0$.
    By Theorem~\ref{thm:inexact_VI}, we obtain
    \begin{align}
        \norm{\hat{P}_i - P^*}_{2} &\le \beta(\norm{\hat{P}_0 - P^*}_{2}, i) + \gamma(\norm{\Delta}_{\infty})
        \notag\\&
        \le \beta(\norm{\hat{P}_0 - P^*}_{2}, i) + \gamma(c(\bar{\mathcal{C}}(\hat{P}_{0}))\norm{\Delta Q_{i}}_{\infty}) \,.
    \end{align}
    To conclude the proof, set the $\mathcal{K}$-function to $\tilde{\gamma}(r) \triangleq \gamma(c(\bar{\mathcal{C}}(\hat{P}_{0})) r)$.

\subsection{Proof of Proposition \ref{prop:robust_LSVI}}\label{apdx:proof_prop_rlsvi}
For all $\hat{P}_{i}\in\mathcal{P}(\hat{P}_{0})$, the sequence is bounded. Indeed by construction, the sequence $\{\Delta Q_{i}\}_{i\geq0}$ induces $\{\Delta_{i}\}_{i\geq0}$ such that $\norm{\Delta}_{\infty}<d^{*}\leq \eigmin{S}$. Hence $\mathcal{P}(\hat{P}_{0})\subset\mathcal{C}(\hat{P}_{0})$, and by Corollary \ref{cor:compactC}, $\bar{\mathcal{C}}(\hat{P}_{0})$ is a compact set. Therefore, $\bar{\mathcal{P}}(\hat{P}_{0})$ is also a compact set.
From Algorithm \ref{alg:R-LSVI} we have that 
\begin{equation}
\label{eq:proof_2.2}
    \hat{Q}_{i} = \sveci{\mathcal{E} (\Theta^{\dagger}_{T}\Psi_{T}\svec{\hat{P}_{i}} + \Theta^{\dagger}_{T}\Xi_{T} )} \,.
\end{equation}
Let $\mathcal{T}$ be the operator defined by
\begin{equation}
    \hat{Q}_{i} = \mathcal{T} (\Theta_{T}, \Psi_{T}, \Xi_{T}, \hat{P}_{i} ) \,.
\end{equation}
Recall that ${Q}_{i} = \mathcal{Q}(\hat{P}_{i})$ and $\Delta Q_{i} = \hat{Q}_{i} - {Q}_{i}$, and hence 
\begin{equation}
    \Delta Q_{i} =  \mathcal{T} (\Theta_{T}, \Psi_{T}, \Xi_{T}, \hat{P}_{i} ) -  \mathcal{T} (\Theta, \Psi, \Xi, \hat{P}_{i} )
\end{equation}
where 
\begin{equation}
    \mathcal{T}(\Theta, \Psi, \Xi, \hat{P}_{i})  =  \left(\begin{array}{cc}
    A^{\top} \hat{P}_{i} A+S & A^{\top} \hat{P}_{i} B \\
    B^{\top} \hat{P}_{i} A & B^{\top} \hat{P}_{i} B+R
    \end{array}\right) \,.
\end{equation}
By the ergodic theorem, we have 
\begin{equation}
\lim _{N \rightarrow \infty} \Theta_{T}=\Theta,\quad
\lim _{N \rightarrow \infty} \Psi_{T}=\Psi,\quad
\lim _{N \rightarrow \infty} \Xi_{T}=\Xi  
\end{equation}
with $\Theta_{T}$ invertible by Assumption \ref{assum:invertible1}. By the continuity of matrix inversion, there exists $\epsilon>0$ such that for all $\Theta_{T}\in\mathcal{B}_{\epsilon}\left(\Theta\right)$, $\Theta_{T}$ is invertible. Hence there exists $T_{1}$ such that for all $T>T_{1}$, $\Theta_{T}\in\mathcal{B}_{\epsilon}\left(\Theta\right)$.  Furthermore, the function $ \mathcal{T} (x, y, z, p)$ is continuous in $(x,y,z)$ for $x$ invertible, and is analytic in $p$, therefore, for any $\hat{P}_{i} \in \bar{\mathcal{P}}(\hat{P}_{0})$ and any $d^{*}>0$, there exists $T_{\hat{P}_{i}}>T_{1}$ such that for all $T>T_{\hat{P}_{i}}$
\begin{align}
\label{eq:local_bound}
    \left\|\Delta Q_{i}\right\|_{2} &= \norm{\mathcal{T} (\Theta_{T}, \Psi_{T}, \Xi_{T}, \hat{P}_{i}) - \mathcal{T} (\Theta, \Psi, \Xi, \hat{P}_{i}) }_{2}
    < d^{*} \,.
\end{align}
Finally, because $\bar{\mathcal{P}}(\hat{P}_{0})$ is compact, there exists $T_0\in\Zset_{+}$ such that for all $T>T_0>T_{1}$, equation \eqref{eq:local_bound} holds.

\section{Dynamic Portfolio Allocation}\label{apdx:GP}
% \subsection{Well-Posedness and Stabilizability}
\noindent We demonstrate that the dynamic portfolio allocation problem, introduced in Section~\ref{sec:portfolio-allocation}, is a well-defined stochastic LQR problem.  
\begin{lem}
\label{lem:S_def_pos}
    For $N,M \in \Zset_{+}$, $\gamma>0$ and $\Sigma \in\Sset^{N}_{++}$, we have
    \begin{equation}
        G = \left[\begin{array}{ccc} \gamma \Sigma & 0 & -I_{N} \\ 0 & I_{M} & 0\\ -I_{N} & 0 & 2\gamma^{-1}\Sigma^{-1} \end{array}\right]\in\Sset^{2N + M}_{++}\,.
    \end{equation}
\end{lem}
\begin{proof}
    By construction $G\in\Sset^{2N + M}$, so it remains to show that the spectrum of $G$ satisfies $\sigma(G)\subset \Rset_{+}^{*}$. Defining $\tilde{\Sigma} \triangleq \gamma\Sigma$, we note that $\tilde{\Sigma},\tilde{\Sigma}^{-1} \in\Sset^{N}_{++}$ are simultaneously diagonalizable. As eigenvalues are invariant under change of basis, without loss of generality we consider $\gamma\Sigma = D$ where $D = \operatorname{diag}(d_{1},\dots,d_{N})$ with $d_{i} = \gamma\sigma_{i}^{2} > 0$ and $\{\sigma_{i}^{2}\}_{i=1}^{N}$ are the strictly positive eigenvalues of $\Sigma$.\\
    To determine the eigenvalues of $G$, we need to solve the system
    \begin{equation}
        G\left[x^{\top}, y^{\top}, z^{\top}\right]^{\top} = \lambda \left[x^{\top}, y^{\top}, z^{\top}\right]^{\top}
    \end{equation}
where $x, z \in\Rset^{N},y\in\Rset^{M}$ and $[x^{\top},y^{\top},z^{\top}]^{\top} \neq 0_{2N+M}$. This is equivalent to solving the system
\begin{equation}
    \begin{cases}
    d_{i}x_{i} - z_{i}&=\lambda x_{i},\quad \forall i\in\{1,\ldots,N\}\\
    y_{j}&=\lambda y_{j},\quad \forall j\in\{1,\ldots,M\}\\
    2d_{i}^{-1}z_{i} - x_{i}&=\lambda z_{i},\quad \forall i\in\{1,\ldots,N\}\,.
\end{cases}
\end{equation}
Since we only require the eigenvalues $\lambda$, it suffices to consider the following cases:\\
\textcircled{1} $\lambda = 0$: Then $0_{2N+M}$ is the only solution. Hence $0\notin  \sigma(G)$.\\
\textcircled{2} $\lambda = 1$: Then $[x^{\top},y^{\top},z^{\top}]^{\top} = [0_{N}^{\top},y^{\top},0_{N}^{\top}]^{\top}$ where $y\neq 0_{M}$. Hence $1\in  \sigma(G)$. \\
\textcircled{3} $\lambda \neq 1$: Then $y = 0_{M}$ and the system reduces to 
\begin{align}
    \begin{cases}
    z_{i}&= (d_{i}-\lambda)x_{i},\quad \forall i\in\{1,\ldots,N\}\\
     x_{i}&=(2d_{i}^{-1}-\lambda)(d_{i}-\lambda)x_{i},\quad \forall i\in\{1,\ldots,N\}\,.
\end{cases}
\end{align}
We are looking for $x\neq0_{N}$ (otherwise this leads to $[x^{\top},y^{\top},z^{\top}]^{\top} = 0_{2N+M}$). Let $i_{0}\in\{1,\cdots,N\}$ be an index such that $x_{i_{0}}\neq 0$. Then $\lambda$ is the positive solution of 
\begin{equation}
    \lambda^{2} - (2 d_{i}^{-1} + d_{i})\lambda + 1 =0
\end{equation}
given by $\lambda = \frac{1}{2} \Big((2 d_{i}^{-1} + d_{i}) \pm \sqrt{(2 d_{i}^{-1})^{2} + (d_{i})^{2}}\Big) > 0$, where the strict positivity follows from  $(a+b)^{2}>a^{2} + b^{2}$ for $a,b >0$. Hence $\sigma(G)\subset \Rset_{+}^{*}$.
\end{proof}
\begin{lem}
    Define $A,B,S$ as in \eqref{eq:matrix_gp1}--\eqref{eq:matrix_gp2} with $\star_{f} = I_{M}$ and $\star_{r} = 2\gamma^{-1}\Sigma^{-1}$, then
    \begin{enumerate}[(i)]
        \item $S\in \Sset^{2N+M}_{++}\,.$
        \item $(A, B)$ is stabilizable.
    \end{enumerate}
\end{lem}
\begin{proof}
\begin{enumerate}[(i)]
     \item This is a direct application of Lemma \ref{lem:S_def_pos}.
     \item  We need to prove that $\exists K \in \Rset^{N \times (2N + M)}$ s.t. $\rho\left(A-B K\right)<1$. Setting $ K = [I_{N}, 0_{N\times M},  0_{N \times N}] $, we obtain  
     \begin{equation}
        A - BK = \left(\begin{array}{ccc}
                    0 & 0 & 0 \\
                    0 & I_{M} - \Phi & 0 \\
                    0 & \Pi & 0
                    \end{array}\right)
    \end{equation}
    and hence $\sigma(A - BK) = \{0\} \bigcup \sigma(I_{M} - \Phi)$, where by assumption $\rho(I_{M} - \Phi) <1$.
\end{enumerate}
\end{proof}

% \endgroup 
% End of reduced spacing
\bibliographystyle{plain}
\bibliography{ref}

@article{kolm2022mean,
	author = {Kolm, Petter N. and Westray, Nicholas},
	journal = {Journal of Portfolio Management},
	number = {6},
	title = {Mean-Variance Optimization for Simulation of Order Flow},
	volume = {48},
	year = {2022}}

@article{garleanu2013dynamic,
	author = {G{\^a}rleanu, N. and Pedersen, L. H.},
	journal = {The Journal of Finance},
	number = {6},
	pages = {2309--2340},
	publisher = {Wiley Online Library},
	title = {Dynamic Trading with Predictable Returns and Transaction Costs},
	volume = {68},
	year = {2013}}

@book{aastrom1995adaptive,
	address = {MA, USA},
	author = {{\AA}str{\"o}m, K. J. and Wittenmark, B.},
	edition = {2nd},
	publisher = {Addison-Wesley},
	title = {Adaptive Control},
	year = {1997}}

@article{kingma2014adam,
  title={Adam: A method for stochastic optimization},
  author={Kingma, Diederik P and Ba, Jimmy},
  journal={arXiv preprint arXiv:1412.6980},
  year={2014}
}

@book{montague1999reinforcement,
	address = {Cambridge, MA},
	author = {Sutton, R. S. and Barto, A. G.},
	edition = {2nd},
	publisher = {MIT Press},
	title = {Reinforcement Learning: An Introduction},
	year = {2018}}

@book{teschl2024ordinary,
	address = {Providence, Rhode Island},
	author = {Teschl, G.},
	publisher = {American Mathematical Society},
	title = {Ordinary Differential Equations and Dynamical Systems},
	volume = {140},
	year = {2024}}

@book{magnus2019matrix,
	address = {Hoboken, New Jersey},
	author = {Magnus, J. R. and Neudecker, H.},
	publisher = {John Wiley \& Sons},
	title = {Matrix Differential Calculus with Applications in Statistics and Econometrics},
	year = {2019}}

@article{beard1997galerkin,
	author = {Beard, R. W. and Saridis, G. N. and Wen, J. T.},
	journal = {Automatica},
	number = {12},
	pages = {2159--2177},
	publisher = {Elsevier},
	title = {Galerkin Approximations of the Generalized {Hamilton-Jacobi-Bellman} Equation},
	volume = {33},
	year = {1997}}

@article{leake1967construction,
	author = {Leake, R. J. and Liu, R.-W.},
	journal = {SIAM Journal on Control},
	number = {1},
	pages = {54--63},
	publisher = {SIAM},
	title = {Construction of Suboptimal Control Sequences},
	volume = {5},
	year = {1967}}

@article{kleinman1969optimal,
	author = {Kleinman, D.},
	journal = {IEEE Transactions on Automatic Control},
	number = {6},
	pages = {673--677},
	publisher = {IEEE},
	title = {Optimal Stationary Control of Linear Systems with Control-Dependent Noise},
	volume = {14},
	year = {1969}}

@article{kleinman1968riccati,
	author = {Kleinman, D.},
	journal = {IEEE Transactions on Automatic Control},
	number = {1},
	pages = {114--115},
	publisher = {IEEE},
	title = {On an Iterative Technique for Riccati Equation Computations},
	volume = {13},
	year = {1968}}

@book{puterman2014markov,
	address = {Hoboken, New Jersey},
	author = {Puterman, M. L.},
	publisher = {John Wiley \& Sons},
	title = {Markov Decision Processes: Discrete Stochastic Dynamic Programming},
	year = {2014}}

@book{bertsekas2011dynamic,
	address = {Belmont, MA},
	author = {Bertsekas, D. P.},
	edition = {3rd},
	publisher = {Athena Scientific},
	title = {Dynamic Programming and Optimal Control},
	volume = {1},
	year = {2011}}

@book{lewis2012optimal,
	address = {Hoboken, New Jersey},
	author = {Lewis, F. L. and Vrabie, D. and Syrmos, V. L.},
	publisher = {John Wiley \& Sons},
	title = {Optimal Control},
	year = {2012}}

@book{stevens2015aircraft,
	address = {Hoboken, New Jersey},
	author = {Stevens, B. L. and Lewis, F. L. and Johnson, E. N.},
	publisher = {John Wiley \& Sons},
	title = {Aircraft Control and Simulation: Dynamics, Controls Design, and Autonomous Systems},
	year = {2015}}

@article{kiumarsi2017optimal,
	author = {Kiumarsi, B. and Vamvoudakis, K. G. and Modares, H. and Lewis, F. L.},
	journal = {IEEE Transactions on Neural Networks and Learning Systems},
	number = {6},
	pages = {2042--2062},
	publisher = {IEEE},
	title = {Optimal and Autonomous Control Using Reinforcement Learning: A Survey},
	volume = {29},
	year = {2017}}

@book{kamalapurkar2018reinforcement,
	address = {Berlin},
	author = {Kamalapurkar, R. and Walters, P. and Rosenfeld, J. and Dixon, W.},
	publisher = {Springer},
	title = {Reinforcement Learning for Optimal Feedback Control},
	year = {2018}}

@article{jiang2012computational,
	author = {Jiang, Y. and Jiang, Z.-P.},
	journal = {Automatica},
	number = {10},
	pages = {2699--2704},
	publisher = {Elsevier},
	title = {Computational Adaptive Optimal Control for Continuous-Time Linear Systems with Completely Unknown Dynamics},
	volume = {48},
	year = {2012}}

@article{jiang2014adaptive,
	author = {Jiang, Y. and Jiang, Z.-P.},
	journal = {Biological Cybernetics},
	number = {4},
	pages = {459--473},
	publisher = {Springer},
	title = {Adaptive Dynamic Programming as a Theory of Sensorimotor Control},
	volume = {108},
	year = {2014}}

@article{jiang2020learning,
	author = {Jiang, Z.-P. and Bian, T. and Gao, W.},
	journal = {Foundations and Trends{\textregistered} in Systems and Control},
	number = {3},
	pages = {176--284},
	publisher = {Now Publishers, Inc.},
	title = {Learning-Based Control: A Tutorial and Some Recent Results},
	volume = {8},
	year = {2020}}

@article{bian2014adaptive,
	author = {Bian, T. and Jiang, Y. and Jiang, Z.-P.},
	journal = {Automatica},
	number = {10},
	pages = {2624--2632},
	publisher = {Elsevier},
	title = {Adaptive Dynamic Programming and Optimal Control of Nonlinear Nonaffine Systems},
	volume = {50},
	year = {2014}}

@article{bian2016adaptive,
	author = {Bian, T. and Jiang, Y. and Jiang, Z.-P.},
	journal = {IEEE Transactions on Automatic Control},
	number = {12},
	pages = {4170--4175},
	publisher = {IEEE},
	title = {Adaptive Dynamic Programming for Stochastic Systems with State and Control Dependent Noise},
	volume = {61},
	year = {2016}}

@article{bian2016value,
	author = {Bian, T. and Jiang, Z.-P.},
	journal = {Automatica},
	pages = {348--360},
	publisher = {Elsevier},
	title = {Value Iteration and Adaptive Dynamic Programming for Data-Driven Adaptive Optimal Control Design},
	volume = {71},
	year = {2016}}

@article{bian2019continuous,
	author = {Bian, T. and Jiang, Z.-P.},
	journal = {SIAM Journal on Control and Optimization},
	number = {6},
	pages = {4150--4174},
	publisher = {SIAM},
	title = {Continuous-Time Robust Dynamic Programming},
	volume = {57},
	year = {2019}}

@article{bian2020model,
	author = {Bian, T. and Wolpert, D. M. and Jiang, Z.-P.},
	journal = {Neural Computation},
	number = {3},
	pages = {562--595},
	publisher = {MIT Press},
	title = {Model-Free Robust Optimal Feedback Mechanisms of Biological Motor Control},
	volume = {32},
	year = {2020}}

@article{Bian2021VI,
	author = {Bian, T. and Jiang, Z.-P.},
	doi = {10.1109/TNNLS.2020.3045087},
	journal = {IEEE Transactions on Neural Networks and Learning Systems},
	number = {7},
	pages = {2781--2790},
	title = {Reinforcement Learning and Adaptive Optimal Control for Continuous-Time Nonlinear Systems: A Value Iteration Approach},
	volume = {33},
	year = {2022}}

@inproceedings{pang2021robust,
	author = {Pang, B. and Jiang, Z. -P.},
	booktitle = {Proceedings of the AAAI conference on artificial intelligence},
	number = {10},
	pages = {9303--9311},
	title = {Robust reinforcement learning: A case study in linear quadratic regulation},
	volume = {35},
	year = {2021}}

@article{pang2022reinforcement,
	author = {Pang, B. and Jiang, Z.-P.},
	journal = {IEEE Transactions on Automatic Control},
	number = {4},
	pages = {2383--2390},
	publisher = {IEEE},
	title = {Reinforcement Learning for Adaptive Optimal Stationary Control of Linear Stochastic Systems},
	volume = {68},
	year = {2022}}

@inproceedings{abbasi2011regret,
	author = {Abbasi-Yadkori, Y. and Szepesv{\'a}ri, C.},
	booktitle = {Proceedings of the 24th Annual Conference on Learning Theory},
	organization = {JMLR Workshop and Conference Proceedings},
	pages = {1--26},
	title = {Regret Bounds for the Adaptive Control of Linear Quadratic Systems},
	year = {2011}}

@inproceedings{abbasi2019model,
	author = {Abbasi-Yadkori, Y. and Lazic, N. and Szepesv{\'a}ri, C.},
	booktitle = {The 22nd International Conference on Artificial Intelligence and Statistics},
	organization = {PMLR},
	pages = {3108--3117},
	title = {Model-Free Linear Quadratic Control via Reduction to Expert Prediction},
	year = {2019}}

@inproceedings{tu2018least,
	author = {Tu, S. and Recht, B.},
	booktitle = {International Conference on Machine Learning},
	organization = {PMLR},
	pages = {5005--5014},
	title = {Least-Squares Temporal Difference Learning for the Linear Quadratic Regulator},
	year = {2018}}

@article{recht2019tour,
	author = {Recht, B.},
	journal = {Annual Review of Control, Robotics, and Autonomous Systems},
	pages = {253--279},
	publisher = {Annual Reviews},
	title = {A Tour of Reinforcement Learning: The View from Continuous Control},
	volume = {2},
	year = {2019}}

@techreport{gao2014machine,
  author      = {Gao, Jim},
  title       = {Machine learning applications for data center optimization},
  institution = {Google},
  address     = {Mountain View, CA, USA},
  type        = {White Paper},
  number      = {21},
  year        = {2014}
}

@article{krauth2019finite,
	author = {Krauth, K. and Tu, S. and Recht, B.},
	journal = {Advances in Neural Information Processing Systems},
	title = {Finite-Time Analysis of Approximate Policy Iteration for the Linear Quadratic Regulator},
	volume = {32},
	year = {2019}}

@inproceedings{tu2019gap,
	author = {Tu, S. and Recht, B.},
	booktitle = {Conference on Learning Theory},
	organization = {PMLR},
	pages = {3036--3083},
	title = {The Gap Between Model-Based and Model-Free Methods on the Linear Quadratic Regulator: An Asymptotic Viewpoint},
	year = {2019}}

@article{dean2020sample,
	author = {Dean, S. and Mania, H. and Matni, N. and Recht, B. and Tu, S.},
	journal = {Foundations of Computational Mathematics},
	number = {4},
	pages = {633--679},
	publisher = {Springer},
	title = {On the Sample Complexity of the Linear Quadratic Regulator},
	volume = {20},
	year = {2020}}

@article{willems2005note,
	author = {Willems, J. C. and Rapisarda, P. and Markovsky, I. and De Moor, B. L. M.},
	journal = {Systems \& Control Letters},
	number = {4},
	pages = {325--329},
	publisher = {Elsevier},
	title = {A Note on Persistency of Excitation},
	volume = {54},
	year = {2005}}

@article{hewer1971iterative,
	author = {Hewer, G.},
	journal = {IEEE Transactions on Automatic Control},
	number = {4},
	pages = {382--384},
	publisher = {IEEE},
	title = {An Iterative Technique for the Computation of the Steady State Gains for the Discrete Optimal Regulator},
	volume = {16},
	year = {1971}}

@inproceedings{ha2023automated,
	author = {Ha, W. Y. and Chakraborty, S. and Yu, Y. and Ghasemi, S. and Jiang, Z.-P.},
	booktitle = {IEEE 26th International Conference on Intelligent Transportation Systems (ITSC)},
	pages = {4215--4220},
	title = {Automated Lane Changing Through Learning-Based Control: An Experimental Study},
	year = {2023}}

@article{jiang1994small,
	author = {Jiang, Z.-P. and Teel, A. R. and Praly, L.},
	journal = {Mathematics of Control, Signals and Systems},
	pages = {95--120},
	publisher = {Springer},
	title = {Small-Gain Theorem for {ISS} Systems and Applications},
	volume = {7},
	year = {1994}}

@article{zhang2023revisiting,
	author = {Zhang, X. and Ba{\c s}ar, T.},
	journal = {IEEE Control Systems Letters},
	pages = {1664-1669},
	publisher = {IEEE},
	title = {Revisiting {LQR} Control from the Perspective of Receding-Horizon Policy Gradient},
	volume = {7},
	year = {2023}}

@article{Lee2022TAC,
	author = {Lee, D.},
	doi = {10.1109/TAC.2022.3181752},
	journal = {IEEE Transactions on Automatic Control},
	number = {10},
	pages = {5661--5668},
	title = {Convergence of Dynamic Programming on the Semidefinite Cone for Discrete-Time Infinite-Horizon {LQR}},
	volume = {67},
	year = {2022}}

@article{Cui2024SCL,
	author = {Cui, L. and Jiang, Z.-P. and Sontag, E. D.},
	doi = {10.1016/j.sysconle.2024.105804},
	journal = {Systems \& Control Letters},
	pages = {105804},
	title = {Small-Disturbance Input-to-State Stability of Perturbed Gradient Flows: Applications to {LQR} Problem},
	volume = {188},
	year = {2024}}

@inproceedings{Cui2023CDC,
	author = {Cui, L. and Pang, B. and Jiang, Z.-P.},
	booktitle = {62nd IEEE Conference on Decision and Control (CDC)},
	doi = {10.1109/CDC49753.2023.10384286},
	pages = {7944--7949},
	title = {Reinforcement-Learning-Based Risk-Sensitive Optimal Feedback Mechanisms of Biological Motor Control},
	year = {2023}}

@article{Cui2021robot,
	author = {Cui, L. and Wang, S. and Zhang, J. and Zhang, D. and Lai, J. and Zheng, Y. and Zhang, Z. and Jiang, Z.-P.},
	doi = {10.1109/LRA.2021.3100269},
	journal = {IEEE Robotics and Automation Letters},
	number = {4},
	pages = {7667--7674},
	title = {Learning-Based Balance Control of Wheel-Legged Robots},
	volume = {6},
	year = {2021}}

@inproceedings{Barto1994,
	author = {Bradtke, S. J. and Ydstie, B. E. and Barto, A. G.},
	booktitle = {Proceedings of 1994 American Control Conference},
	doi = {10.1109/ACC.1994.735224},
	pages = {3475--3479},
	title = {Adaptive Linear Quadratic Control Using Policy Iteration},
	year = {1994}}

@book{howard1960dynamic,
	address = {New York},
	author = {Howard, R. A.},
	publisher = {John Wiley \& Sons},
	title = {Dynamic Programming and {M}arkov Processes},
	year = {1960}}

@book{bellman1966dynamic,
	address = {Princeton, NJ},
	author = {Bellman, R.},
	publisher = {Princeton University Press},
	title = {Dynamic Programming},
	year = {1957}}

@article{Huang2022,
	author = {Huang, M. and Jiang, Z.-P. and Ozbay, K.},
	doi = {10.1109/TCYB.2020.3029077},
	journal = {IEEE Transactions on Cybernetics},
	number = {6},
	pages = {5267--5277},
	title = {Learning-Based Adaptive Optimal Control for Connected Vehicles in Mixed Traffic: Robustness to Driver Reaction Time},
	volume = {52},
	year = {2022}}

@article{Kalman1960,
	author = {Kalman, R.},
	journal = {Boletín de la Sociedad Matemática Mexicana},
	number = {2},
	pages = {102--119},
	title = {Contribution to the Theory of Optimal Control},
	volume = {5},
	year = {1960}}

@inbook{Sontag2008,
	address = {Germany},
	author = {Sontag, E. D.},
	booktitle = {Nonlinear and Optimal Control Theory},
	pages = {163--220},
	publisher = {Springer Verlag},
	series = {Lecture Notes in Mathematics},
	title = {Input-to-State Stability: Basic Concepts and Results},
	year = {2008}}

@article{Sontag1989,
	author = {Sontag, E. D.},
	journal = {IEEE Transactions on Automatic Control},
	number = {4},
	pages = {435--443},
	title = {Smooth Stabilization Implies Coprime Factorization},
	volume = {34},
	year = {1989}}

@book{book_Jiang,
	address = {Hoboken, New Jersey},
	author = {Jiang, Y. and Jiang, Z.-P.},
	publisher = {Wiley-IEEE Press},
	title = {Robust Adaptive Dynamic Programming},
	year = {2017}}

@article{song2024convergence,
  title={Convergence and robustness of value and policy iteration for the linear quadratic regulator},
  author={Song, Bowen and Wu, Chenxuan and Iannelli, Andrea},
  journal={arXiv preprint arXiv:2411.04548},
  year={2024}
}

@article{abeille2016lqg,
  title={LQG for portfolio optimization},
  author={Abeille, Marc and Lazaric, Alessandro and Brokmann, Xavier and others},
  journal={arXiv preprint arXiv:1611.00997},
  year={2016}
}

@article{moallemi2017dynamic,
  title={Dynamic portfolio choice with linear rebalancing rules},
  author={Moallemi, Ciamac C and Sa{\u{g}}lam, Mehmet},
  journal={Journal of Financial and Quantitative Analysis},
  volume={52},
  number={3},
  pages={1247--1278},
  year={2017},
  publisher={Cambridge University Press}
}

@article{lai2025value,
  title={Value Iteration for Stochastic {LQR} With Convergence Guarantees},
  author={Lai, Jing and Xiong, Junlin and Kang, Yu},
  journal={IEEE Transactions on Neural Networks and Learning Systems},
  year={2025},
  publisher={IEEE}
}

@article{fan2024value,
  title={Value iteration for {LQR} control of unknown stochastic-parameter linear systems},
  author={Fan, Wenwu and Xiong, Junlin},
  journal={Systems \& Control Letters},
  volume={185},
  pages={105731},
  year={2024},
  publisher={Elsevier}
}

@article{shen2024data,
  title={Data-driven near optimization for fast sampling singularly perturbed systems},
  author={Shen, Hao and Peng, Chuanjun and Yan, Huaicheng and Xu, Shengyuan},
  journal={IEEE Transactions on Automatic Control},
  volume={69},
  number={7},
  pages={4689--4694},
  year={2024},
  publisher={IEEE}
}

@article{yi2021ALQG,		
	title   = {Adaptive Linear Quadratic Control for Stochastic Discrete-Time Linear Systems with Unmeasurable Multiplicative and Additive Noises},
	author  = {Jiang, Yi and Liu, Lu and Feng, Gang},
	journal = {IEEE Transactions on Automatic Control},
	year    = {2024},   
	volume  = {69},    
	number  = {11},    
	pages   = {7808-7815},
	month   = {November},
	doi     = {10.1109/TAC.2024.3399630},
}

@article{yi2021reinforcementlearning,		
	title   = {Adaptive Optimal Control of Networked Control Systems with Two-Channel Stochastic Dropouts},
	author  = {Jiang, Yi and Gao, Weinan and Chen, Ci and Chai, Tianyou and Lewis, Frank L},
	journal = {SIAM Journal on Control and Optimization},
	year    = {2023},   
	volume  = {61},    
	number  = {5},    
	pages   = {3183-3208},
	month   = {October},
	doi     = {10.1137/21M1438797},
}

@book{horn2012matrix,
  title={Matrix analysis},
  author={Horn, Roger A and Johnson, Charles R},
  year={2012},
  publisher={Cambridge University Press}
}
\end{document}